\documentclass[12pt]{amsart}
\usepackage{macros}
\usepackage{mathtools}

\usepackage{todonotes}

\newcommand{\gmdn}[1]{{g^{\{\underline{#1}\}^d} }}
\newcommand{\dec}[1]{\downarrow_{#1}\!}

\title[]{Tableau posets and the fake degrees of coinvariant algebras}

\author{Sara C. Billey, Matja\v{z} Konvalinka, Joshua P. Swanson}
\address{Billey: Department of Mathematics, University of Washington,
  Seattle, WA 98195, USA}
\email{billey@math.washington.edu}

\thanks{The first author was partially supported by the Washington
  Research Foundation and  DMS-1764012. The second author was partially supported by Research Project
BI-US/16-17-042 of the Slovenian Research Agency and research core funding No. P1-0294.}
\address{Konvalinka: Faculty of Mathematics and Physics,
University of Ljubljana,
Jadranska 21, Ljubljana, Slovenia, and Institute for Mathematics, Physics and Mechanics, Jadranska 19, Ljubljana, Slovenia}
\email{matjaz.konvalinka@fmf.uni-lj.si}
\address{Swanson: Department of Mathematics,
University of California, San Diego (UCSD),
La Jolla, CA  92093-0112}
\email{jswanson@ucsd.edu}
\date{\today}

\begin{document}

\begin{abstract}
  We introduce two new partial orders on the standard Young tableaux
  of a given partition shape, in analogy with the strong and weak
  Bruhat orders on permutations.  Both posets are ranked by the major
  index statistic offset by a fixed shift.  The existence of such
  ranked poset structures allows us to classify the realizable major
  index statistics on standard tableaux of arbitrary straight shape
  and certain skew shapes.  By a theorem of Lusztig--Stanley, this
  classification can be interpreted as determining which irreducible
  representations of the symmetric group exist in which homogeneous
  components of the corresponding coinvariant algebra, strengthening a
  recent result of the third author for the modular major index.  Our
  approach is to identify patterns in standard tableaux that allow one
  to mutate descent sets in a controlled manner.  By work of Lusztig
  and Stembridge, the arguments extend to a classification of all
  nonzero fake degrees of coinvariant algebras for finite complex
  reflection groups in the infinite family of Shephard--Todd groups.
\end{abstract}
\maketitle

\section{Introduction}
\label{sec:intro}

Let $\SYT(\lambda)$ denote the set of all standard Young tableaux of
partition shape $\lambda$.  We say $i$ is a \emph{descent} in a
standard tableau $T$ if $i+1$ comes before $i$ in the row reading word
of $T$, read from bottom to top along rows in English notation.
Equivalently, $i$ is a descent in $T$ if $i+1$ appears in a lower row
in $T$.  Let $\maj(T)$ denote the \emph{major index statistic} on
$\SYT(\lambda)$, which is defined to be the sum of the descents of
$T$.  The major index generating function for $\SYT(\lambda)$ is given
by
\begin{equation}\label{eq:syt-q}
    \SYT(\lambda)^{\maj}(q) := \sum_{T \in \SYT(\lambda)} q^{\maj(T)} =
    \sum_{k\geq 0} b_{\lambda,k} q^k.
\end{equation}

The polynomial $\SYT(\lambda)^{\maj}(q)$ has two elegant closed forms,
one due to Steinberg based on dimensions of
irreducible representations of $\GL_n(\mathbb{F}_q)$, see \cite{steinberg.1951}, and one due to
Stanley \cite {stanley.1979} generalizing the Hook-Length Formula; see
\Cref{thm:stanley_maj}.

For fixed $\lambda$, consider the \emph{fake degree sequence}
\begin{equation}\label{eq:a_la}
  b_{\lambda, k} := \#\{T \in \SYT(\lambda) : \maj(T) = k\} \text{ for }
  k=0,1,2,\ldots
\end{equation}
The fake degrees have appeared in a variety of algebraic and
representation-theoretic contexts  including Green's work on the
degree polynomials of uni\-po\-tent $\GL_n(\bF_q)$-representations
\cite[Lemma~7.4]{MR0072878}, the irreducible de\-com\-position of type $A$
coinvariant algebras \cite[Prop.~4.11]{stanley.1979}, Lusztig's work
on the irreducible representations of classical groups
\cite{Lusztig.77}, and branching rules between symmetric groups and
cyclic subgroups \cite[Thm.~3.3]{stembridge89}.  The term ``fake
degree'' was apparently coined by Lusztig \cite{Carter.89}, perhaps
because $\#\SYT(\lambda)= \sum_{k\geq 0} b_{\lambda,k}$ is the degree
of the irreducible $S_n$-representation indexed by $\lambda$, so a
$q$-analog of this number is not itself a degree but related to the
degree.

We consider three natural enumerative questions involving the fake degrees:
\begin{enumerate}[(I)]
  \item which $b_{\lambda,k}$ are zero?
  \item are the fake degree sequences unimodal?
  \item are there efficient asymptotic estimates for $b_{\lambda, k}$?
\end{enumerate}

We completely settle (I) with the following result. Denote by $\lambda'$
the con\-ju\-gate partition of $\lambda$, and let $b(\lambda)
:= \sum_{i \geq 1} (i-1) \lambda_i$.

\begin{Theorem}\label{thm:zeros}
  For every partition $\lambda\vdash n\geq1$ and integer $k$ such that
  $b(\lambda) \leq k \leq \binom{n}{2}-b(\lambda')$, we have
  $b_{\lambda, k} > 0$ except in the case when $\lambda$ is
  a rectangle with at least two rows and columns and $k$ is either
  $b(\lambda)+1$ or $\binom{n}{2}-b(\lambda')-1$. Furthermore,
  $b_{\lambda, k} = 0$ for $k < b(\lambda)$ or $k > \binom{n}{2}
  - b(\lambda')$.
\end{Theorem}

As a consequence of the proof of \Cref{thm:zeros}, we identify two ranked
poset structures on $\SYT(\lambda)$ where the rank function is
determined by $\maj$.  Furthermore, as a corollary of \Cref{thm:zeros}
we have a new proof of a complete classification due to the third author
\cite[Thm.~1.4]{s17} generalizing an earlier result of Klyachko
\cite{klyachko74} for when the counts
\[ a_{\lambda, r} := \#\{T \in \SYT(\lambda) : \maj(T) \equiv_n r\} \]
for $\lambda \vdash n$ are nonzero.

The easy answer to question (II) is ``no''.  The fake degree sequences
are not always unimodal.  For example, $ \SYT(4,2)^{\maj}(q)$ is not
unimodal.  See \Cref{Ex:42}. Nonetheless, certain inversion
number generating functions $p_\alpha^{(k)}(q)$ which appear in
a generalization of $\SYT(\lambda)^{\maj}(q)$ are in
fact unimodal; see \Cref{def:pam} and \Cref{cor:pam_unimodal}.
Furthermore, computational evidence suggests $\SYT(\lambda)^{\maj}(q)$
is typically not far from unimodal.

Questions (II) and (III) are addressed in a separate article
\cite{bks2}.  In particular, we show in that article that
the coefficients of $\SYT(\lambda^{(i)})^{\maj}(q)$ are asymptotically
normal for any sequence of partitions
$\lambda^{(1)}, \lambda^{(2)}, \ldots$ such that $\aft(\lambda^{(i)})$
approaches infinity where $\aft(\lambda)$ is the number of boxes
outside the first row or column, whichever is smaller.  The aft
statistics on partitions is in FindStat as \cite[St001214]{findstat}.

We note that there are polynomial expressions for the fake degrees
$b_{\lambda,k}$ in terms of parameters $H_i$, the number of cells of
$\lambda$ with hook length equal to $i$.  These polynomials are
closely related to polynomials that express the number of permutations
$S_n$ of a given inversion number $k\leq n$ as a function of $n$ by
work of Knuth. See \Cref{poly-h} and \Cref{cor:mahonian.numbers}.
These polynomials are useful in some cases, however, we find that in
practice Stanley's formula is the most effective way to compute a
given fake degree sequence for partitions up to size 200.  See
\Cref{rem:cyclotomic} for more on efficient computation using
cyclotomic polynomials.

Symmetric groups are the finite reflection groups of type $A$.  The
clas\-si\-fi\-cation and invariant theory of both finite irreducible real
reflection groups and complex reflection groups developed over the
past century builds on our understanding of the type $A$ case
\cite{Hum}.  In particular, these groups are classified by
Shephard--Todd into an infinite family $G(m,d,n)$ together with $34$ exceptions.
Using work of Stembridge on generalized exponents for irreducible
representations, the analog of \eqref{eq:syt-q}
can be phrased for all Shephard--Todd groups as
\begin{equation}\label{eq:gmdn}
    \gmdn{\lambda}(q) \coloneqq
      \frac{\#\{\underline{\lambda}\}^d}{d} \cdot
      \bbinom{n}{\alpha(\underline{\lambda})}_{q; d} \cdot
      \prod_{i=1}^m \SYT(\lambda^{(i)})^{\maj}(q^m) =
      \sum b_{\{\underline{\lambda}\}^d, k} q^k
    \end{equation}
    where
    $\underline{\lambda}=(\lambda^{(1)}, \ldots, \lambda^{(m)})$ is
    a sequence of $m$ partitions with $n$ cells total,
    $\alpha(\underline{\lambda})=(|\lambda^{(1)}|,
    \ldots,|\lambda^{(m)}|)\vDash n$, $d \mid m$, and
    $\{\underline{\lambda}\}^d$ is the orbit of $\underline{\lambda}$
    under the group $C_d$ of $(m/d)$-fold cyclic rotations;
    see \Cref{lem:g_la_d}.
The polynomials
    $ \bbinom{n}{\alpha(\underline{\lambda})}_{q; d} $ are
    deformations of the usual $q$-multinomial coefficients which we
    explore in \Cref{sec:deformed}.
    The
    coefficients $b_{\{\underline{\lambda}\}^d, k}$ are the fake
    degrees in this case.

    We use \eqref{eq:gmdn} and \Cref{thm:zeros} to completely classify
    all nonzero fake degrees for coinvariant algebras for all
    Shephard--Todd groups $G(m,d,n)$, which in\-cludes the finite
    real reflection groups in types $A$, $B$, and $D$.  See
    \Cref{cor:internal_zeros_B} and \Cref{cor:internal_zeros_D} for the
    type $B$ and $D$ cases, respectively.  See
    \Cref{thm:block_internal_zeros} and \Cref{thm:block_internal_zeroes.2}
    for the general $C_m \wr S_n$ and $G(m, d, n)$ cases, respectively.

    The rest of the paper is organized as follows. In
    \Cref{sec:background}, we give back\-ground on tableau
    combinatorics, Shephard--Todd groups, and their irreducible
    representations.  \Cref{sec:h-poly} describes the polynomial
    formulas for fake degrees in type $A$.  \Cref{sec:internal_zeros}
    presents our combinatorial argument proving \Cref{thm:zeros} and
    giving poset structures on tableaux of a given shape. \Cref{sec:AER}
    uses the argument in \Cref{sec:internal_zeros} to answer in the
    affirmative a question of Adin--Elizalde--Roichman about internal
    zeros of $\SYT(\lambda)^{\des}(q)$; see \Cref{cor:AER}. In
    \Cref{ssec:zeros_wreath}, we begin to address the question of
    characterizing nonzero fake degrees by starting with the wreath
    products $C_m \wr S_n=G(m,1,n)$; see The\-o\-rem \ref{thm:block_internal_zeros}.
    In \Cref{sec:deformed}, we define the deformed $q$-multinomials
    $\bbinom{n}{\alpha}_{q; d}$ as rational functions and give a
    summation formula, \Cref{thm:q_d_mult_sum}, which shows they are
    polynomial.  Finally, in \Cref{sec:other_types}, we complete the
    classification of nonzero fake degrees for $G(m,d,n)$ and spell
    out how \eqref{eq:gmdn} relates to Stembridge's original
    generating function for the fake degrees in $G(m,d,n)$; see
    \Cref{thm:block_internal_zeroes.2} and \Cref{lem:g_la_d}.
    We discuss potential algebraic and geometric directions
    for future work in \Cref{sec:future}.

\section{Background}\label{sec:background}

In this section, we review some standard terminology and results on
com\-bi\-na\-to\-ri\-al statistics and tableaux.  Many further details
in this area can be found in \cite{ec1,ec2}.  We also review
background on the finite complex reflection groups and their
irreducible representations.  Further details in this area can be
found in \cite{Carter.89,Sagan.1991}.

\subsection{Word and Tableau Combinatorics}\label{ssec:words_tableaux}

Here we review standard com\-bi\-na\-to\-ri\-al notions related to
words and tableaux.

\begin{Definition}
  Given a word $w = w_1w_2\cdots w_n$ with letters
  $w_i \in \bZ_{\geq 1}$, the \textit{content} of $w$ is the sequence
  $\alpha = (\alpha_1, \alpha_2, \ldots)$ where $\alpha_i$ is the
  number of times $i$ appears in $w$.  Such a sequence $\alpha$ is
  called a (weak) \textit{composition} of $n$, written as
  $\alpha \vDash n$. Trailing $0$'s are often omitted when writing
  compositions, so $\alpha = (\alpha_1, \alpha_2, \ldots, \alpha_m)$
  for some $m$.  Note, a word of content $(1,1,\ldots,1)\vDash n$ is a
  permutation in the symmetric group $S_n$ written in one-line
  notation. The \textit{inversion number} of $w$ is
    \[ \inv(w) := \#\{(i, j) : i < j, w_i > w_j\}. \]
  The \textit{descent set} of $w$ is
    \[ \Des(w) := \{0<i <n : w_i > w_{i+1}\} \]
  and the \textit{major index} of $w$ is
    \[ \maj(w) := \sum_{i \in \Des(w)} i. \]
\end{Definition}

The study of permutation statistics is a classical topic in enumerative
com\-bi\-na\-tor\-ics.  The major index statistic on
permutations was introduced by Percy MacMahon in his seminal works
\cite{MacMahon.1913,MR1576566}.  At first glance, this function on
permutations may be unintuitive, but it has inspired hundreds of
papers and many generalizations; for example on Macdonald polynomials
\cite{HHL2005}, posets \cite{Ehrenborg.Readdy.2015}, quasisymmetric
functions \cite{Shareshian-Wachs.2010}, cyclic sieving
\cite{Reiner-Stanton-White.CSP,Ahlbach-Swanson.2017}, and bijective
combinatorics \cite{Foata,carlitz.1975}.

\begin{Definition}
  Given a finite set $W$ and a function
  $\stat \colon W \to \bZ_{\geq 0}$, write the corresponding ordinary
  generating function as
    \[ W^{\stat}(q) := \sum_{w \in W} q^{\stat(w)}. \]
\end{Definition}

\begin{Definition}
  Let $\alpha = (\alpha_1, \ldots, \alpha_m) \vDash n$. We use the
  following standard $q$-analogues:
\[
  \begin{array}{lllll}  \vspace{.2in}
    [n]_q &:=& 1 + q + \cdots + q^{n-1} = \frac{q^n - 1}{q - 1}, &
                                                                  \hspace{.3in}
    & \text{($q$-integer)}\\ \vspace{.2in}
    {}
   [n]_q! &:=& [n]_q [n-1]_q \cdots [1]_q, &
                                                                  \hspace{.3in}
    & \text{($q$-factorial)}\\ \vspace{.2in}

    \binom{n}{k}_q
      &:=& \frac{[n]_q!}{[k]_q! [n-k]_q!} \in \bZ_{\geq 0}[q], &
                                                                  \hspace{.3in}
    & \text{($q$-binomial)}\\ \vspace{.2in}
    \binom{n}{\alpha}_q
      &:=& \frac{[n]_q!}{[\alpha_1]_q! \cdots [\alpha_m]_q!} \in
           \bZ_{\geq 0}[q]
           &
                                                                  \hspace{.3in}
    & \text{($q$-multinomial).}\\
  \end{array}
  \]
\end{Definition}

\begin{Example}
  The identity statistic on the set $W = \{0, \ldots, n-1\}$ has
  gen\-er\-ating function $[n]_q$. The ``sum'' statistic on
  $W = \prod_{j=1}^n \{0, \ldots, j-1\}$ has generating function
  $[n]_q!$. It is straightforward to show that also
  $S_n^{\inv}:=\sum_{w\in S_n}q^{\inv(w)}= [n]_q!$.
\end{Example}

For $\alpha \vDash n$, let $\W_\alpha$ denote the set of all words of content
$\alpha$.  A classic result of MacMahon is that $\maj$ and $\inv$ have
the same distribution on $\W_\alpha$ which is determined by the
corresponding $q$-multinomial.

\begin{Theorem}{\cite[\S1]{MR1576566}}\label{thm:macmahon}
  For each $\alpha \vDash n$,
  \begin{align}\label{eqn:macmahon}
     \W_\alpha^{\maj}(q) = \binom{n}{\alpha}_q = \W_\alpha^{\inv}(q).
  \end{align}
\end{Theorem}

\begin{Definition}
  A polynomial $P(q) = \sum_{i=0}^n c_i q^i$ of degree $n$ is \emph{symmetric} if
  $c_i = c_{n-i}$ for $0\leq i \leq n$.  We generally say $P(q)$ is
  \emph{symmetric} also if there exists an integer $k$ such that
  $q^k P(q)$ is symmetric.  We say $P(q)$ is \emph{unimodal} if
  \[ c_0\leq c_1 \leq \cdots \leq c_j \geq c_{j+1} \geq \cdots \geq
    c_n \] for some $0\leq j \leq n$.  Furthermore, $P(q)$ has
  \emph{no internal zeros} provided that $c_j \neq 0$ whenever
  $c_i, c_k \neq 0$ and $i<j<k$.

\end{Definition}

From \Cref{thm:macmahon} and the definition of the $q$-multinomials,
we see that each $\W_\alpha^{\maj}(q)$ is a symmetric polynomial with
constant and leading coefficient $1$.  Indeed, these polynomials are
\textit{unimodal} generalizing the well-known case for Gaussian
coefficients \cite[Thm 3.1]{Stanley.1980} and \cite{Zeilberger.1986}.
It also follows easily from MacMahon's theorem that $\W_\alpha^{\maj}(q)$ has
no internal zeros.

\subsection{Partitions and Standard Young Tableaux}

\begin{Definition}
  A composition $\lambda \vDash n$ such that
  $\lambda_1 \geq \lambda_2 \geq \ldots$ is called a
  \textit{partition} of $n$, written as $\lambda \vdash n$. The
  \textit{size} of $\lambda$ is $|\lambda| := n$ and the
  \textit{length} $\ell(\lambda)$ of $\lambda$ is the number of
  non-zero entries. The \textit{Young diagram} of $\lambda$ is the
  upper-left justified arrangement of unit squares called
  \textit{cells} where the $i$th row from the top has $\lambda_i$
  cells following the English notation; see \Cref{fig:partition_a}.
  The cells of a tableau are indexed by matrix notation when we refer
  to their row and column.  The \textit{hook length} of a cell
  $c \in \lambda$ is the number $h_c$ of cells in $\lambda$ in the
  same row as $c$ to the right of $c$ and in the same column as $c$
  and below $c$, including $c$ itself; see \Cref{fig:partition_c}.  A
  \emph{corner} of $\lambda$ is any cell with hook length $1$.
  A \emph{notch} of $\lambda$ is any $(i,j)$ not in $\lambda$
  such that both $(i-1,j)$ and $(i,j-1)$ are in $\lambda$. Note that
  notches cannot be in the first row or column of $\lambda$.  A
  \emph{bijective filling} of $\lambda$ is any labeling of the cells
  of $\lambda$ by the numbers $[n]=\{1,2,\ldots,n\}$.  The symmetric
  group $S_n$ acts on bijective fillings of $\lambda$ by acting on the
  labels.
\end{Definition}

\begin{figure}[ht]
  \centering
  \begin{subfigure}[t]{0.32\textwidth}
    \centering
    \[
    \tableau{{} & {} & {} & {} & {} & {} \\ {} & {} & {} \\ {} & {} & {} }
    \]
    \caption{Young diagram of $\lambda$.}
    \label{fig:partition_a}
  \end{subfigure}
  \hspace{.2cm}
  \begin{subfigure}[t]{0.32\textwidth}
    \centering
    \[
    \tableau{{8} & {7} & {6} & {3} & {2} & {1} \\ {4} & {3} & {2} \\ {3} & {2} & {1} }
    \]
    \caption{Hook lengths of $\lambda$.}
    \label{fig:partition_c}
  \end{subfigure}
  \caption{Constructions related to the partition $\lambda=(6, 3, 3)\vdash 12$.
    The partition has corners at positions $(3, 3)$ and $(1, 6)$ and one notch at position $(2, 4)$.}
  \label{fig:partition}
\end{figure}

\begin{Definition}
  A \textit{skew partition} $\lambda/\nu$ is a pair of partitions
  $(\nu, \lambda)$ such that the Young diagram of $\nu$ is contained
  in the Young diagram of $\lambda$. The cells of $\lambda/\nu$ are
  the cells in the diagram of $\lambda$ which are not in the diagram
  of $\nu$, written $c \in \lambda/\nu$.  We identify straight
  partitions $\lambda$ with skew partitions $\lambda/\varnothing$
  where $\varnothing = (0, 0, \ldots)$ is the empty partition. The
  \textit{size} of $\lambda/\nu$ is $|\lambda/\nu| := |\lambda| - |\nu|$. The
  notions of bijective filling, hook lengths, corners, and notches
  naturally extend to skew partitions as well.
\end{Definition}

\begin{figure}[ht]
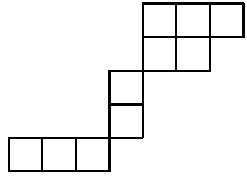

  \centering
    \[
    \tableau{& & & & {} & {} & {} \\ & & & & {} & {} \\ & & & {} \\ & & & {} \\ {} & {} & {} }
    \]
    \caption{Diagram for the skew partition
      $\lambda/\nu = 76443/4433$, which is also the block diagonal
      skew shape
      $\underline{\lambda} = ((3, 2), (1, 1),(3) )$.}
    \label{fig:diag_a}
\end{figure}

\begin{Definition}\label{def:block.diag.skew.partions}
  Given a sequence of partitions $\underline{\lambda} = (\lambda^{(1)},
  \ldots, \lambda^{(m)})$, we identify the sequence with the
  \emph{block diagonal skew partition} obtained by translating the Young diagrams of the
  $\lambda^{(i)}$ so that the rows and columns occupied by these
  components are disjoint, form a valid skew shape, and they appear in order from
  top to bottom as depicted in  \Cref{fig:diag_a}.
\end{Definition}

\begin{figure}[ht]
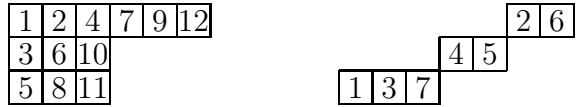

  \centering
  \begin{subfigure}[t]{0.32\textwidth}
    \centering
    \[
    \tableau{{1} & {2} & {4} & {7} & {9} & {12} \\ {3} & {6} & {10} \\ {5} & {8} & {11} }
    \]
  \end{subfigure}
    \begin{subfigure}[t]{0.32\textwidth}
    \centering
    \[
    \tableau{& & & & & 2 & 6 \\ & & & 4 & 5 \\ 1 & 3 & 7}
    \]
  \end{subfigure}
  \caption{On the left is a standard Young tableau of straight shape
    $\lambda=(6,3,3)$ with descent set  $\{2, 4, 7, 9, 10\}$ and major
    index $32$.  On the right is a standard Young tableau of block diagonal
    skew shape $(7,5,3)/(5,3)$ corresponding to the sequence of partitions
    $((2), (2), (3))$ with descent set $\{2,6\}$ and major index $8$. }
  \label{fig:SYTpartition}
\end{figure}

  \begin{Definition}
    A \textit{standard Young tableau} of shape $\lambda/\nu$ is a
    bijective filling of the cells of $\lambda/\nu$ such that labels
    increase to the right in rows and down columns; see
    \Cref{fig:SYTpartition}.   The set of standard Young
    tableaux of shape $\lambda/\nu$ is denoted $\SYT(\lambda/\nu)$.
    The \textit{descent set} of $T \in \SYT(\lambda/\nu)$ is the set
    $\Des(T)$ of all labels $i$ in $T$ such that $i+1$ is in a
    strictly lower row than $i$.  The \textit{major index} of $T$ is
    \[ \maj(T) := \sum_{i \in \Des(T)} i. \]
\end{Definition}

The block diagonal skew partitions $\underline{\lambda}$ allow us to
simultaneously consider words and tableaux as follows.  Let
$\W_\alpha$ be the set of all words with content $\alpha
= (\alpha_1, \ldots, \alpha_k)$.  Letting
$\underline{\lambda} = ((\alpha_k), \ldots, (\alpha_1))$, we have a bijection
\begin{equation}\label{eq:bij}
  \phi \colon \SYT(\underline{\lambda}) \too{\sim} \W_\alpha
\end{equation}
which sends a tableau $T$ to the word whose $i$th letter
is the row number in which $i$ appears in $T$, counting from the \textit{bottom up}
rather than top down.
For example, using the skew tableau $T$ on the right of \Cref{fig:SYTpartition}, we have
$\phi(T) = 1312231 \in \W_{(3, 2, 2)}$. It is easy to see that
$\Des(\phi(T)) = \Des(T)$, so that $\maj(\phi(T))
= \maj(T)$.

\subsection{Major Index Generating Functions}\label{ssec:maj_gfs}


 Stanley gave the following an\-a\-logue of
\Cref{thm:macmahon} for standard Young tableaux of a given shape. It
generalizes the famous Frame--Robinson--Thrall Hook-Length Formula
\cite[Thm. 1]{Frame-Robinson-Thrall.1954} or \cite[Cor.~7.21.6]{ec2} obtained by setting $q=1$.

\begin{Theorem}{\cite[7.21.5]{ec2}}\label{thm:stanley_maj}
  Let $\lambda \vdash n$ with $\lambda = (\lambda_1, \lambda_2, \ldots)$. Then
  \begin{equation}\label{eq:stanley_maj}
    \SYT(\lambda)^{\maj}(q)
       = \frac{q^{b(\lambda)}[n]_q!}
                  {\prod_{c \in \lambda} [h_c]_q}
  \end{equation}
  where $b(\lambda) := \sum (i-1)\lambda_i$ and $h_c$ is the hook
  length of the cell $c$.
\end{Theorem}

\begin{Remark}\label{rem:cyclotomic}
  Since $\#\SYT(\lambda)$ typically grows extremely quickly, Stanley's
  formula offers a practical way to compute
  $\SYT(\lambda)^{\maj}(q)$ even when $n \approx 100$
  by expressing both the numerator and
  denominator, up to a $q$-shift, as a product of cyclotomic
  polynomials and canceling all factors from the denominator.  We prefer
  to use cyclotomic factors over linear factors in order to avoid arithmetic
  in cyclotomic fields.
\end{Remark}

\begin{Example}\label{Ex:42}
  For $\lambda=(4,2)$, $b(\lambda)=2$ and the multiset of hook lengths
  is $\{1^2,2^2,4,5\}$ so $|\SYT(\lambda)| = 9$ by the Hook-Length
  Formula.  The major index generating function is given by
  \begin{align*}
    \SYT(4,2)^{\maj}(q) &=   q^8 + q^7 + 2q^6 + q^5 + 2q^4 + q^3 + q^2\\
                        &= q^2 \frac{[6]_q!}{[5]_q[4]_q[2]_q[2]_q} =
q^2 \frac{[6]_q[3]_q}{[2]_q}.
  \end{align*}
  Note, $\SYT(4,2)^{\maj}(q)$ is symmetric but not unimodal.

  For $\lambda=(4,2,1)$, $b(\lambda)=4$ and the multiset of hook
  lengths is $\{1^3, 2,3,4,6\}$ so $|\SYT(\lambda)| =35$ by the Hook-Length
  Formula.  The major index generating function is given by
  \begin{align*}
  \SYT(4,2,1)^{\maj}(q) &= q^{14} + 2q^{13} + 3q^{12} + 4q^{11} +
                            5q^{10} + 5q^{9} + 5q^{8} + 4q^{7} + 3q^{6} \\
                            &+ 2q^{5} + q^{4} = q^4 \frac{[7]_q!}{[6]_q[4]_q[3]_q[2]_q}=q^4[7]_q[5]_q.
  \end{align*}
  Note, $\SYT(4,2,1)^{\maj}(q)$ is symmetric and unimodal.
\end{Example}

\begin{Example}
  We recover $q$-integers, $q$-binomials, and $q$-Catalan numbers,
  up to $q$-shifts as special cases of the major index generating
  function for tableaux as follows:
  \begin{align*}
    \SYT(\lambda)^{\maj}(q)
      &= \begin{cases}
             q[n]_q & \text{if }\lambda = (n, 1), \\
             q^{\binom{k+1}{2}} \binom{n}{k}_q & \text{if }\lambda = (n-k+1, 1^k), \\
             q^n \frac{1}{[n+1]_q} \binom{2n}{n}_q & \text{if }\lambda = (n, n).
           \end{cases}
  \end{align*}
\end{Example}

The following strengthening of Stanley's formula to $\underline{\lambda}$
is well known (e.g.\ see \cite[(5.6)]{stembridge89}), though since it is
somewhat difficult to find explicitly in the literature, we include a short proof.

\begin{Theorem}\label{thm:diag_maj} 
  Let $\underline{\lambda} = (\lambda^{(1)}, \ldots, \lambda^{(m)})$
  where $\lambda^{(i)} \vdash \alpha_i$ and $n = \alpha_1 + \cdots + \alpha_m$.
  Then
  \begin{equation}\label{eq:diag_maj}
    \SYT(\underline{\lambda})^{\maj}(q)
      = \binom{n}{\alpha_1, \ldots, \alpha_m}_q \cdot \prod_{i=1}^m \SYT(\lambda^{(i)})^{\maj}(q).
  \end{equation}

  \begin{proof}
    The stable principal specialization of skew Schur
    functions is given by
      \[ s_{\lambda/\nu}(1, q, q^2, \ldots) = \frac{\SYT(\lambda/\nu)^{\maj}(q)}
          {\prod_{j=1}^{|\lambda/\nu|} (1 - q^j)}; \]
    see \cite[Lemma~3.1]{stembridge89} or \cite[Prop.7.19.11]{ec2}. On the other hand,
    it is easy to see from the definition of a skew Schur function as
    the content generating function for semistandard tableaux of the
    given shape  that
      \[ s_{\underline{\lambda}}(x_1, x_2, \ldots)
          = \prod_{i=1}^m s_{\lambda^{(i)}}(x_1, x_2, \ldots). \]
    The result quickly follows.
  \end{proof}
\end{Theorem}

\begin{Remark}
  \Cref{thm:stanley_maj} and \Cref{thm:diag_maj} have several immediate
  cor\-ol\-laries.  First, we recover MacMahon's result,
  \Cref{thm:macmahon}, from \Cref{thm:diag_maj} when
  $\underline{\lambda} = ((\alpha_m), (\alpha_{m-1}), \ldots)$ by using the
  $\maj$-preserving bijection $\phi$ in \eqref{eq:bij}. Second, each
  $\SYT(\underline{\lambda})^{\maj}(q)$ is symmetric (up to a $q$-shift)
  with leading co\-ef\-ficient $1$. In particular, there is a unique
  ``maj-minimizer'' and ``maj-max\-i\-mizer'' tableau in each
  $\SYT(\underline{\lambda})$. Moreover,
  \begin{equation}\label{eq:min.maj}
    \min \maj(\SYT(\underline{\lambda}))
    = b(\underline{\lambda})
    \end{equation}
  and
  \begin{equation}\label{eq:max.maj}
    \max \maj(\SYT(\underline{\lambda})) = \binom{n}{2} - b(\underline{\lambda}')
      = b(\underline{\lambda}) + \binom{|\underline{\lambda}|+1}{2}
      - \sum_{c \in \underline{\lambda}} h_c
    \end{equation}
    where $b(\underline{\lambda}) := \sum_i b(\lambda^{(i)})$ and
    $b(\underline{\lambda}') := \sum_i b(\lambda^{(i) \, \prime})$.
\end{Remark}

       For general skew shapes, $\SYT(\lambda/\nu)^{\maj}(q)$ does not
       factor as a product of cyclotomic polynomials times $q$ to a
       power. A ``$q$-Naruse'' formula due to Morales--Pak--Panova,
       \cite[(3.4)]{1512.08348}, gives an analogue of
       \Cref{thm:stanley_maj} involving a sum over ``excited
       diagrams,'' though the resulting sum has a single term
       precisely for the block diagonal skew partitions $\underline{\lambda}$.

\subsection{Complex Reflection Groups}\label{sub:complex.reflection.groups}

A \textit{complex reflection group} is a finite subgroup of
$\GL(\bC^n)$ generated by \textit{pseudo-reflections}, which are
elements which pointwise fix a codimension-$1$
hyperplane. Shephard--Todd, building on work of Coxeter and others, famously
classified the complex reflection groups \cite{Shephard-Todd.54}. The
irreducible representations were constructed by Young, Specht,
Lusztig, and
others. We now summarize these results and fix some notation.

\begin{Definition}
  A \textit{pseudo-permutation matrix} is a matrix where each row and
  column has a single non-zero entry.  For positive integers $m,n$,
  the \emph{wreath product} $C_m \wr S_n \subset \GL(\bC^n)$ is the
  group of $n \times n$ pseudo-permutation matrices whose non-zero
  entries are complex $m$th roots of unity. For $d \mid m$, let
  $G(m, d, n)$ be the \emph{Shephard--Todd group} consisting of
  matrices $x \in C_m \wr S_n$ where the product of the non-zero
  entries in $x$ is an $(m/d)$th root of unity.  In fact, $G(m, d, n)$
  is a normal subgroup of $C_m \wr S_n$ of index $d$ with cyclic
  quotient $(C_m \wr S_n)/G(m, d, n) \cong C_d$ of order $d$.
\end{Definition}

\begin{Theorem}\cite{Shephard-Todd.54}
  Up to isomorphism, the complex reflection groups are precisely the
  direct products of the groups $G(m, d, n)$, along with $34$
  exceptional groups.
\end{Theorem}

\begin{Remark}
  Special cases of the Shephard--Todd groups include the fol\-low\-ing.
  The Weyl group of type $A_{n-1}$, or equivalently the symmetric
  group $S_n$, is isomorphic to $G(1, 1, n)$.  The Weyl groups of both
  types $B_n$ and $C_n$ are $G(2, 1, n)$, the group of $n\times n$
  signed permutation matrices. The subgroup of the group of signed
  permutations whose elements have evenly many negative signs is the
  Weyl group of type $D_n$, or $G(2, 2, n)$ as a Shephard--Todd
  group. We also have that $G(m, m, 2)$ is the dihedral group of order
  $2m$, and $G(m, 1 ,1)$ is the cyclic group $C_m$ of order $m$.
\end{Remark}

The complex irreducible representations of $S_n$ were constructed by
Young \cite{Young.collected.work} and are well known to be certain
modules $S^\lambda$ \textit{canonically} indexed by partitions
$\lambda \vdash n$. These representations are beautifully described in
\cite{Sagan.1991}.  Specht extended the construction to irreducibles
for $G \wr S_n$ where $G$ is a finite group.

\begin{Theorem}\cite{Specht.1935}\label{thm:specht}
  The complex inequivalent irreducible representations of
  $C_m \wr S_n$ are certain modules $S^{\underline{\lambda}}$ indexed
  by the sequences of partitions
  $\underline{\lambda} = (\lambda^{(1)}, \ldots, \lambda^{(m)})$
  for which $|\underline{\lambda}|
  \coloneqq |\lambda^{(1)}| + \cdots + |\lambda^{(m)}| = n$.
\end{Theorem}

\begin{Remark}
  The version we give of \Cref{thm:specht} was stated by
  Stembridge \cite[Thm. 4.1]{stembridge89}.  The $C_m$-irreducibles
  are naturally though non-canonically indexed by $\bZ/m$ up to one of
  $\phi(m)$ additive automorphisms, where $\phi(m)$ is Euler's totient
  function. Correspondingly, one may identify $\bZ/m$ with
  $\{1, \ldots, m\}$ and obtain $\phi(m)$ different indexing schemes
  for the $C_m \wr S_n$-ir\-re\-duc\-ibles. The resulting indexing schemes
  are rearrangements of one another, and our results will be
  independent of these choices.
\end{Remark}

Clifford described a method for determining the branching rules of
ir\-re\-duc\-ible representations for a normal subgroup of a given finite
group \cite{Clifford.1937}.  Stembridge combined this method with
Specht's theorem to describe the ir\-re\-duc\-ible representations for all
Shephard--Todd groups from the $C_m \wr S_n$-irreducible
representations.

We use Stembridge's terminology where possible.  In particular, for
$d \mid m$, the $(m/d)$-\emph{fold cyclic rotations} are the elements
in the subgroup isomorphic to $C_d$ of $S_m$ generated by
$\sigma_m^{m/d}$, where $\sigma_m = (1,2,\ldots,m)$ is the long cycle.  Let
$S_m$ act on block diagonal partitions of the form
$\underline{\lambda} = (\lambda^{(1)}, \ldots, \lambda^{(m)})$ by
permuting the blocks.  This action restricts to
$C_d=\langle \sigma_m^{m/d} \rangle$ as well.  Let
$\{\underline{\lambda}\}^d$ denote the orbit of $\underline{\lambda}$
under the $(m/d)$-fold cyclic rotations in $C_d$.  Note, the number of
block diagonal partitions in such a $C_d$-orbit, denoted
$\#\{\underline{\lambda}\}^d$, always divides $d$, but could be less
than $d$ if $\underline{\lambda}$ contains repeated partitions.

For example, take $d=2$ and $m=6$.  If
$\underline{\lambda} =((1),(2),(3,2),(4),(5),(6,1))$, then
$\{\underline{\lambda}\}^2$ has two elements, $\underline{\lambda}$ and
$((4),(5),(6,1),(1),(2),(3,2))$.  If
$\underline{\mu} =((1),(2),(3,2),(1),(2),(3,2))$, then
$\{\underline{\mu}\}^2$ only contains the element $\underline{\mu}$.

\begin{Theorem}\cite[Remark after Prop.~6.1]{stembridge89}\label{thm:stembridge.irreps}
  The complex inequivalent irreducible representations of $G(m, d, n)$
  are certain modules $S^{\{\underline{\lambda}\}^d, c}$ indexed by
  the pairs $(\{\underline{\lambda}\}^d, c)$ where
  $\underline{\lambda} = (\lambda^{(1)}, \ldots, \lambda^{(m)})$ is
  a sequence of partitions with $|\underline{\lambda}| = n$,
  $\{\underline{\lambda}\}^d$ is the orbit of $\underline{\lambda}$
  under $(m/d)$-fold cyclic rotations, and $c$ is any positive integer
  $1\leq c \leq \frac{d}{\#\{\underline{\lambda}\}^d}.$
\end{Theorem}

\begin{Remark}
  As with $C_m \wr
  S_n$, the indexing scheme is again non-canonical in general up to a
  choice of orbit representative, though our results relying on this
  work are independent of these choices. In fact, Stembridge
  uses $\underline{\lambda} = (\lambda^{(m-1)}, \ldots, \lambda^{(0)})$,
  which is the most natural setting for \Cref{thm:stembridge.irreps}
  and \Cref{thm:canonical_tableaux} below.
  The fake degrees for irreducibles $S^{\underline{\lambda}}$ of $C_m \wr S_n$
  are invariant up to a $q$-shift under all permutations of $\underline{\lambda}$ in $S_m$, so
  for our purposes the indexing scheme is largely irrelevant. The fake degrees
  for irreducibles $S^{\{\underline{\lambda}\}^d, c}$ of $G(m, d, n)$, however,
  are only invariant under the $(m/d)$-fold cyclic rotations of $\underline{\lambda}$
  in general. In this case, strictly
  speaking our $\lambda^{(i)}$
  corresponds to the irreducible cyclic group representation $\chi^{i-1}$ defined by
  $\chi^{i-1}(\sigma_m) = \omega_m^{i-1}$ where $\omega_m$ is a fixed
  primitive $m$th root of unity in the sense that
    \[ S^{\underline{\lambda}} \coloneqq \left(\chi^0 \wr S^{\lambda^{(1)}}
        \otimes \cdots \otimes \chi^{m-1}
        \wr S^{\lambda^{(m)}}\right)
        \ind_{C_m \wr S_{\alpha(\underline{\lambda})}}^{C_m \wr S_n}; \]
  see \cite[(4.1)]{stembridge89}. Since we have no need of these
  explicit representations, we
  have used the naive indexing scheme throughout.
\end{Remark}

\begin{Example}
  For the type $B_n$ group $G(2, 1, n)$, the irreducible
  rep\-re\-sen\-ta\-tions are indexed by pairs $(\lambda, \mu)$ since $C_1$ is
  the trivial group and so in each case $c=1$.
\end{Example}

\begin{Example}
  For the type $D_n$ group $G(2, 2, n)$, the irreducible
  rep\-re\-sen\-ta\-tions can be thought of as being indexed by the sets
  $\{\lambda, \mu\}$ with $\lambda \neq \mu$ and
  $|\lambda| + |\mu| = n$ together with the pairs $(\nu, 1)$ and
  $(\nu, 2)$ where $\nu \vdash n/2$. The orbits alone can be thought
  of as the 2 element multisets $\{\lambda, \mu\}$ with
  $|\lambda| + |\mu| = n$.
\end{Example}

\subsection{Coinvariant Algebras}
As mentioned in the introduction, Stanley (see \cite{stanley.1979}) and
Lusztig (unpublished) determined the graded irreducible de\-com\-pos\-ition
of the type $A$ coinvariant algebra via the major index generating
function on standard Young tableaux. Stembridge was the first to
publish a complete proof of this result and extended it to the complex
reflection groups $G(m, d, n)$ \cite{stembridge89}. We now summarize
these results.

\begin{Definition}
  Any group $G \subset \GL(\bC^n)$ acts on the polynomial ring with
  $n$ variables $\bC[x_1, \ldots, x_n]$ by identifying $\bC^n$ with
  $\Span_{\bC}\{x_1, \ldots, x_n\}$ and extending the $G$-action
  multiplicatively.  The \textit{coinvariant algebra} of $G$ is the
  quotient of $\bC[x_1, \ldots, x_n]$ by the ideal generated by
  homogeneous $G$-invariant polynomials of positive degree, which is
  thus a graded $G$-module.
\end{Definition}

\begin{Definition}
  Let $R_n$ denote the coinvariant algebra of $S_n$. For $\lambda \vdash n$,
  let $g^\lambda(q)$ be the \textit{fake degree polynomial} whose $k$th coefficient
  is the multiplicity of $S^\lambda$ in the $k$th degree piece of $R_n$.
\end{Definition}

\begin{Theorem}[Lusztig--Stanley, {\cite[Prop.~4.11]{stanley.1979}}]
  For a partition $\lambda$,
    \[ g^\lambda(q) = \SYT(\lambda)^{\maj}(q). \]
  Equivalently, the multiplicity of $S^\lambda$ in the $k$th degree piece of the type $A$
  coinvariant algebra $R_n$ is $b_{\lambda, k}$, the number of standard tableaux of
  shape $\lambda \vdash n$ with major index $k$.
\end{Theorem}

\begin{Definition}
  Let $R_{m, n}$ denote the coinvariant algebra of $C_m \wr S_n$.
  Set
    \[ b_{\underline{\lambda}, k} \coloneqq
        \text{the multiplicity of $S^{\underline{\lambda}}$ in the
                $k$th degree piece of $R_{m, n}$}. \]
  Write the corresponding \emph{fake degree polynomial} as
    \[ g^{\underline{\lambda}}(q) \coloneqq \sum_k b_{\underline{\lambda}, k} q^k. \]
\end{Definition}

\begin{Definition}
  Given a sequence of partitions $\underline{\lambda} =
  (\lambda^{(1)}, \ldots, \lambda^{(m)})$, recall
  \[ b(\alpha(\underline{\lambda})) = \sum_{i=1}^m (i-1)
    |\lambda^{(i)}|. \] We continue to identify $\underline{\lambda}$
  with a block diagonal skew partition when con\-ve\-nient, as in
  \Cref{def:block.diag.skew.partions}.  Thus,
  $\SYT(\underline{\lambda})$ is the set of standard Young tableaux on
  the block diagonal skew partition $\underline{\lambda}$.  We will
  abuse notation and define $b(\alpha(T)) \coloneqq b(\alpha(\underline{\lambda}))$ for any
  $T \in \SYT(\underline{\lambda})$, which is not necessary in the
  next theorem but will be essential for the general Shephard--Todd
  groups $G(m,d,n)$.
\end{Definition}

\begin{Theorem}\cite[Thm.~5.3]{stembridge89}\label{thm:block_exponents}
  For $\underline{\lambda} = (\lambda^{(1)}, \ldots, \lambda^{(m)})$
  with $|\underline{\lambda}| = n$,
    \[ g^{\underline{\lambda}}(q) = q^{b(\alpha(\underline{\lambda}))}
        \SYT(\underline{\lambda})^{\maj}(q^m). \]
  Equivalently, the multiplicity of $S^{\underline{\lambda}}$ in the $k$th
  degree piece of the $C_m \wr S_n$ coinvariant algebra $R_{m, n}$ is the
  number of standard tableaux $T$ of block diagonal shape $\underline{\lambda}$
  with $k = b(\alpha(T)) + m\cdot\maj(T)$.
\end{Theorem}

\begin{Remark}
  By \eqref{eq:diag_maj}, we have an
  explicit product formula for $g^{\underline{\lambda}}(q)$ also.
  Furthermore, in \cite{bks2}, we characterize the possible
  limiting distributions for the coefficients of the polynomials
  $ \SYT(\underline{\lambda})^{\maj}(q).$ We show that in most
  cases, the limiting distribution is the normal distribution.
  Consequently, that characterization can be interpreted as a
  statement about the asymptotic distribution of irreducible
  components in different degrees of the $C_m \wr S_n$ coinvariant
  algebras.
\end{Remark}

\begin{Corollary}
  In type $B_n$, the irreducible indexed by $(\lambda, \mu)$ with
  $|\lambda| = k$ and $|\mu| = n-k$ has fake degree polynomial
    \[ g^{(\lambda, \mu)}(q)
        = q^{|\mu| + 2b(\lambda) + 2b(\mu)}
           \binom{n}{k}_{q^2}\
           \frac{[k]_{q^2}!}{\prod_{c \in \lambda} [h_c]_{q^2}}
           \frac{[n-k]_{q^2}!}{\prod_{c' \in \mu} [h_{c'}]_{q^2}}. \]
\end{Corollary}

\begin{Definition}
  Let $R_{m, d, n}$ denote the coinvariant algebra of $G(m, d, n)$
  assuming $d \mid m$. For an orbit $\{\underline{\lambda}\}^d$ of a
  sequence of $m$ partitions with $n$ total cells under $(m/d)$-fold cyclic rotations, set
  \[ b_{\{\underline{\lambda}\}^d, k} \coloneqq \text{the
      multiplicity of $S^{\{\underline{\lambda}\}^d, c}$ in the $k$th
      degree piece of $R_{m, d, n}$}, \] which in fact depends only on
  the orbit $\{\underline{\lambda}\}^d$ and not the number $c$ by
  \cite[Prop.~6.3]{stembridge89}. Write the corresponding \emph{fake degree
  polynomial} as
    \[ \gmdn{\lambda}(q) \coloneqq
      \sum_k b_{\{\underline{\lambda}\}^d, k} q^k. \]
\end{Definition}

\begin{Theorem}\cite[Cor.~6.4]{stembridge89}\label{thm:block_exponents.2}
  Let $\{\underline{\lambda}\}^d$ be the orbit of a sequence of
  $m$ partitions $\underline{\lambda}$ with $|\underline{\lambda}| = n$
  under $(m/d)$-fold cyclic rotations. Then
    \[ \gmdn{\lambda}(q)
        = \frac{\left(\{\underline{\lambda}\}^d\right)^{b\circ \alpha}(q)}{[d]_{q^{nm/d}}}
        \SYT(\underline{\lambda})^{\maj}(q^m) \]
  where
    \[ \left(\{\underline{\lambda}\}^d\right)^{b\circ \alpha}(q) \coloneqq
        \sum_{\underline{\mu} \in \{\underline{\lambda}\}^d} q^{b(\alpha(\underline{\mu}))}. \]
\end{Theorem}

\begin{Corollary}[{\cite[Sect. 2.5]{Lusztig.77}}, {\cite[Cor. 6.5]{stembridge89}}]
  In type $D_n$, an irreducible indexed by $\{\underline{\lambda}\}^2$
  with $\underline{\lambda}=(\lambda, \mu)$ and $|\lambda| = k$, $|\mu| = n-k$ has
  fake degree polynomial
    \[ g^{\{\underline{\lambda}\}^2}(q) =
        \kappa_{\lambda\mu} q^{2b(\lambda) + 2b(\mu)}
        \frac{q^k + q^{n-k}}{1+q^n}
        \binom{n}{k}_{q^2}
        \frac{[k]_{q^2}!}{\prod_{c \in \lambda} [h_c]_{q^2}}
        \frac{[n-k]_{q^2}!}{\prod_{c' \in \mu} [h_{c'}]_{q^2}}, \]
  where $\kappa_{\lambda\mu} = 1$ if $\lambda \neq \mu$ and
  $\kappa_{\lambda\lambda} = 1/2$.
\end{Corollary}

Observe that \Cref{thm:block_exponents} gives a direct tableau
interpretation of the coefficients of
$g^{\underline{\lambda}}(q)$. More generally, Stembridge gave a
tableau interpretation of the coefficients of
$\gmdn{\lambda}(q)$ which we next describe.

\begin{Definition}\label{def:cononical.tabs}
  For a given $m,d,n$, let
  $\underline{\lambda}=(\lambda^{(1)},\ldots, \lambda^{(m)}) $ be a
  sequence of $m$ partitions with $|\underline{\lambda}|=n$.  Let
  $\{\underline{\lambda}\}^d$ be the orbit of $\underline{\lambda}$
  under $(m/d)$-fold rotations. The cyclic group
  $C_d = \langle \sigma_m^{m/d} \rangle$ acts on the disjoint union
  $\bigsqcup_{\underline{\mu} \in \{\underline{\lambda}\}^d}
  \SYT(\underline{\mu})$ as follows.  Given
  $\underline{\mu}=(\mu^{(1)}, \ldots, \mu^{(m)}) \in
  \{\underline{\lambda}\}^d$, each $T \in \SYT(\underline{\mu})$ may be
  considered as a sequence
  $\underline{T}=(T^{(1)}, \ldots, T^{(m)})$ of fillings of the
  shapes $\mu^{(i)}$. The group $C_d$ acts by $(m/d)$-fold rotations
  of this sequence of fillings.  Write the resulting orbit as
  $\{\underline{T}\}^d$, which necessarily has size $d$.  For such a
  $\underline{T}$, the largest entry of $\underline{T}$, namely $n$,
  appears in some $T^{(k)}$. If among the elements of the orbit
  $\{\underline{T}\}^d$ of $\underline{T}$ this value $k$ is minimal for
  $\underline{T}$ itself, then we call $\underline{T}$ the
  \textit{canonical standard tableau representative} for
  $\{\underline{T}\}^d$.  Let
  \[\SYT(\{\underline{\lambda}\}^d) \subset
    \bigsqcup_{\underline{\mu} \in \{\underline{\lambda}\}^d}
    \SYT(\underline{\mu})
  \]
  be the \emph{set of canonical standard tableau representatives of
    orbits} $\{\underline{T}\}^d$ for $T \in \SYT(\underline{\lambda})$.  Recall,
  $b(\alpha(\underline{T})):=b(\alpha(\underline{\mu})) = \sum (i-1) |\mu^{(i)}|$ if
    $\underline{T} \in \SYT(\underline{\mu})$, so $b\circ \alpha$ is not generally
    constant on $\SYT(\{\underline{\lambda}\}^d)$.
  \end{Definition}

\begin{Remark}
  When the parts $\lambda^{(i)}$ are all non-empty, the set
  $\SYT(\{\underline{\lambda}\}^d)$ is the set of standard block
  diagonal skew tableaux of some shape
  $\underline{\mu} \in \{\underline{\lambda}\}^d$ where
  $n=|\underline{\lambda}|$ is in the upper-right-most partition
  possible among the $(m/d)$-fold cyclic rotations of its blocks.
  Since every orbit $\{\underline{T}  \}^d$ has size $d$, we have
    \[ \#\SYT(\{\underline{\lambda}\}^d)
        = \frac{\#\{\underline{\lambda}\}^d}{d} \#\SYT(\underline{\lambda}). \]
\end{Remark}

\begin{Theorem}\cite[Thm. 6.6]{stembridge89}\label{thm:canonical_tableaux}
  Let $\underline{\lambda}$ be a sequence of $m$ partitions with
  $|\underline{\lambda}| = n$. Let
  $\{\underline{\lambda}\}^d$ be the orbit of $\underline{\lambda}$ under
  $(m/d)$-fold cyclic rotations. Then
    \[ \gmdn{\lambda}(q) =
        \SYT(\{\underline{\lambda}\}^d)^{b\circ \alpha +m \cdot \maj}(q). \]
  Equivalently, the multiplicity of $S^{\{\underline{\lambda}\}^d, c}$ in the
  $k$th degree piece of the $G(m, d, n)$ coinvariant algebra $R_{m, d, n}$ is
  the number of canonical standard tableaux $\underline{T} \in \SYT(\{\underline{\lambda}\}^d)$
  with $k = b(\alpha(\underline{T})) + m\cdot\maj(\underline{T})$.
\end{Theorem}

\section{Polynomial Formulas For Fake Degrees}\label{sec:h-poly}

In this section, we briefly show how to construct polynomial formulas
for the fake degrees $b_{\lambda, k}$ directly from Stanley's $q$-hook
length formula.  We will use these polynomials in the next section for
small changes from the minimal major index.  Our results extend to a
formula for counting permutations of a given inversion number.

Given $\lambda$, let
\begin{align}
  H_i(\lambda) &= \# \{c \in \lambda \ \colon \ h_c = i\},\\
  m_i (\lambda) &=\#\{k \ \colon \ \lambda_k=i\}.
\end{align}
If $\lambda$ is understood, we abbreviate $H_i=H_i(\lambda)$.  For any
nonnegative integer $k$ and polynomial $f(q)$, let $[q^k]f(q)$ be the
coefficient of $q^k$ in $f(q)$.

\begin{Lemma}\label{poly-h} For every $\lambda \vdash n$ and $k=b(\lambda)+d$, we have
\begin{equation}\label{eq:polys}
b_{\lambda, k}
= [q^{b(\lambda)+d}] \SYT(\lambda)^{\maj}(q) = \sum_{\substack{\mu
    \vdash d \\ \mu_1 \leq n}} \prod_{i=1}^{|\lambda|}
\binom{H_i+m_i(\mu)-2}{m_i(\mu)}
\end{equation}
which is a polynomial in the $H_i$'s for every positive
integer $n$.

\begin{proof}
By \Cref{thm:stanley_maj}, we have
\begin{equation}\label{eq:remark.polys}
q^{-b(\lambda)} \SYT(\lambda)^{\maj}(q)
       = \frac{[n]_q!}
                  {\prod_{c \in \lambda} [h_c]_q} = \prod_{i=1}^n
                  (1-q^i)^{-(H_i-1)}.
\end{equation}
The result follows using the expansion
$(1-q^i)^{-j} = \sum_{n=0}^\infty \binom{j+n-1}{n}q^{in}$ and
mul\-ti\-pli\-ca\-tion of ordinary generating functions.
\end{proof}
\end{Lemma}

  Note that if $H_i(\lambda) = 0$ and $m_i(\mu) = 1$, then the
  corresponding binomial co\-ef\-fi\-cient in \eqref{eq:polys} is $-1$, so it
  is not obvious from this formula that the coefficients
  $b_{\lambda,k}$ are all nonnegative, which is clearly true by
  definition.
\begin{Remark}\label{rem:H_formulas}
  The first few polynomials are given by
  \begin{align*}
    [q^{b(\lambda)+1}] \SYT(\lambda)^{\maj}(q)
      &= H_1 - 1 \\
      &= \#\{\text{notches of $\lambda$}\}, \\
    [q^{b(\lambda)+2}] \SYT(\lambda)^{\maj}(q)
      &= \binom{H_1}{2} + H_2 - 1, \\
    [q^{b(\lambda)+3}] \SYT(\lambda)^{\maj}(q)
      &= \binom{H_1+1}{3} + (H_1 - 1)(H_2 - 1) + (H_3 - 1) \\
    [q^{b(\lambda)+4}] \SYT(\lambda)^{\maj}(q)
      &= \binom{H_1+2}{4} + \binom{H_2}{2} + \binom{H_1}{2} (H_2 - 1) \\
      &\qquad+ (H_1 - 1)(H_3 - 1) + (H_4 - 1).
  \end{align*}
  These exact formulas hold for all $|\lambda| \geq 4$. For smaller
  size partitions some terms will not appear.
\end{Remark}

  It is interesting to compare these polynomials to the ones described
  by Knuth for the number of permutations with $k\leq n$ inversions in
  $S_n$ in \cite[p.16]{Knuth}.  See also \cite[Ex.~1.124]{ec1} and \cite[A008302]{oeis}.  We
  can extend Knuth's formulas to all $0\leq k \leq \binom{n}{2}$ using
  the same idea.

  \begin{Corollary}\label{cor:mahonian.numbers}  For fixed positive integers $k$ and $n$, we have
  \begin{equation}\label{eq:fkn}
     \#\{w\in S_n : \inv(w) =d\} = \sum (-1)^{\#\{\mu_i>1\}}
     \binom{n+ m_1(\mu) - 2}{m_1(\mu)}
    \end{equation}
    where the sum is over all partitions $\mu \vdash d$ such that
    $\mu_1\leq n$ and all of the parts of $\mu$ larger than 1
    are distinct.
\end{Corollary}

The proof follows in exactly the same way from the formula
  \[ \sum_{w \in S_n} q^{\inv(w)}
     = \prod_{i=1}^n[i]_q = \prod_{i=1}^n[i]_q /[1]_q
     = (1-q)^{-n} \prod_{i=1}^n (1-q^i). \]
In essence, this is the case of the
$q$-hook length formula when all of the hooks are of length 1.

\begin{Remark}
  Let $T(d,n)$ be the number of partitions $\mu \vdash d$ such that
  $\mu_1\leq n$ and all of the parts of $\mu$ larger than 1 are
  distinct.  The triangle of numbers $T(d,n)$ for $1\leq n \leq d$ is
  \cite[A318806]{oeis}.
\end{Remark}

\section{Type $A$ Internal Zeros Classification}\label{sec:internal_zeros}
As a corollary of Stanley's formula, we know that for every partition
$\lambda\vdash n\geq 1$ there is a unique tableau with minimal major
index $b(\lambda)$ and a unique tableau with maximal major index
$\binom{n}{2}-b(\lambda')$.  These two agree for shapes consisting of
one row or one column, and otherwise they are distinct.  It is easy to
identify these two tableaux in $\SYT(\lambda)$; see
\Cref{def.min.max.tabs} below.  Then, we classify all of the values
$k$ such that $b(\lambda) < k <\binom{n}{2}-b(\lambda')$ and the fake
degree $b_{\lambda, k}=0$.  We refer to such $k$ as \emph{internal
  zeros}, meaning the location of zeros in the fake degree sequence
for $\lambda$ between the known minimal and maximal nonzero
locations.

\begin{Definition}\label{def.min.max.tabs}
  \
  \begin{enumerate}
  \item The \emph{max-maj tableau} for $\lambda$ is obtained by filling
    the outermost, max\-i\-mum length, vertical strip in $\lambda$ with the
    largest possible numbers
    $|\lambda|,|\lambda|-1,\ldots,|\lambda| - \ell(\lambda) + 1$
    starting from the bottom row and going up, then filling the
    rightmost maximum length vertical strip containing cells not
    previously used with the largest remaining numbers, etc.

  \item The \emph{min-maj tableau} of $\lambda$ is obtained similarly by
    filling the outermost, maximum length,  horizontal strip in
    $\lambda$ with the largest possible num\-bers
    $|\lambda|,|\lambda|-1,\ldots,|\lambda| - \lambda_1+1$ going right
    to left, then filling the lowest maximum length horizontal strip
    containing cells not previously used with the largest remaining
    numbers, etc.
  \end{enumerate}
\end{Definition}

\begin{figure}[ht]
  \centering
  \begin{subfigure}[t]{0.4\textwidth}
    \centering
    \includegraphics{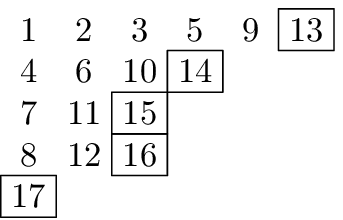}
    \caption{A max-maj tableau and its outermost vertical strip.}
    \label{fig:maxmaj_ex}
  \end{subfigure}
  \hspace{0.05\textwidth}
  \begin{subfigure}[t]{0.4\textwidth}
    \centering
    \includegraphics{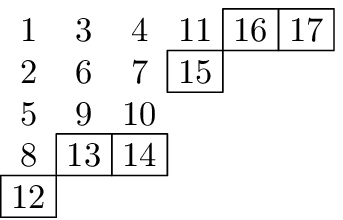}
    \caption{A min-maj tableau and its outermost horizontal strip.}
    \label{fig:minmaj_ex}
  \end{subfigure}
  \caption{Max-maj tableau and min-maj tableau for $\lambda=(6, 4, 3, 3, 1)$.}
  \label{fig:extmaj_ex}
\end{figure}

\noindent
See \Cref{fig:extmaj_ex} for an example. Note that the max-maj
tableau of $\lambda$ is the transpose of the min-maj
tableau of $\lambda'$.

The $q^{b(\lambda)+1}$ coefficients of $\SYT(\lambda)^{\maj}(q)$
can be computed as in \Cref{poly-h} or \Cref{rem:H_formulas}, resulting in the following.

\begin{Corollary}\label{cor:first_few_coeffs}
  We have $[q^{b(\lambda)+1}] \SYT(\lambda)^{\maj}(q) = 0$ if
  and only if $\lambda$ is a rectangle. If $\lambda$ is a rectangle with
  more than one row and column, then
  $[q^{b(\lambda)+2}] \SYT(\lambda)^{\maj}(q) = 1$.
\end{Corollary}

A similar statement holds for $\maj(T)=\binom{n}{2}-b(\lambda')-1$ by
symmetry.  Thus, $\SYT^{\maj}(q)$ has internal zeros when $\lambda$ is
a rectangle with at least two rows and columns.  We will show these are
the only internal zeros of type $A$ fake degrees, proving
\Cref{thm:zeros}.

\begin{Definition}
  Let $\cE(\lambda)$ denote the set of \emph{exceptional}
  tableaux of shape $\lambda$ consisting of the following elements.
  \begin{enumerate}[(i)]
    \item For all $\lambda$, the max-maj tableau for $\lambda$.

    \item If $\lambda$ is a rectangle, the min-maj tableau for
      $\lambda$.

    \item If $\lambda$ is a rectangle with at least two rows
      and columns, the unique tableau in $\SYT(\lambda)$ with
      major index equal to $\binom{n}{2}-b(\lambda')-2$.
      It is obtained from the max-maj tableau of $\lambda$ by
      applying the cycle $(2,3,\ldots,\ell(\lambda)+1)$,
      which reduces the major index by $2$.
  \end{enumerate}
\end{Definition}

For example, $\cE(64331)$ consists of just the max-maj
tableau for $64331$ in \Cref{fig:maxmaj_ex}, while $\cE(555)$ has
the following three elements:

\[
  \begin{matrix}
    1 & 2 & 3 & 4 & 5 \\
    6 & 7 & 8 & 9 & 10 \\
    11 & 12 & 13  & 14 & 15
  \end{matrix}
  \qquad  \qquad
  \begin{matrix}
    1 & 2 & 7 & 10 & 13 \\
    3 & 5 & 8 & 11 &14 \\
    4 & 6 & 9 & 12 & 15
  \end{matrix}
  \qquad \qquad
  \begin{matrix}
    1 & 4 & 7 & 10 & 13 \\
    2 & 5 & 8 & 11 & 14 \\
    3 & 6 & 9 & 12 & 15 \\
  \end{matrix}.
\]

\medskip

We prove \Cref{thm:zeros} by constructing a map
\begin{equation}\label{eq:varphi_telegraph}
  \varphi \, \colon \SYT(\lambda) \setminus \cE(\lambda) \longrightarrow \SYT(\lambda)
\end{equation}
with the property
\begin{equation}\label{eq:varphi_maj}
  \maj(\varphi(T)) = \maj(T) + 1.
\end{equation}
For \emph{most} tableaux $T$, we can find another tableau $T'$ of the
same shape such that $\maj(T')=\maj(T)+1$ by applying some
\emph{simple cycle} to the values in $T$, meaning a permutation whose
cycle notation is either $(i,i+1,\ldots,k-1,k)$ or
$(k,k-1,\ldots, i+1,i)$ for some $i<k$. We will show there are 5
additional rules that must be added to complete the definition.

We note that technically the symmetric group $S_n$ does not act on
$\SYT(\lambda)$ for $\lambda \vdash n$ since this action will not
generally preserve the row and column strict requirements for standard
tableaux.  However, $S_n$ acts on the set of all bijective fillings of $\lambda$
using the alphabet $\{1,2,\ldots, n\}$ by acting on the values.  We
will only apply permutations to tableaux after locating all values in
some interval $[i,j]=\{i,i+1,\ldots, j\}$ in $T$.  The reader is
encouraged to verify that the specified permutations always maintain
the row and column strict properties.

\subsection{Rotation Rules}

We next describe certain configurations in a tableau which imply that
a simple cycle will increase $\maj$ by $1$.  Recall, the cells of
a tableau are indexed by matrix notation.

\begin{Definition}\label{def:pos_rot}
  Given $\lambda \vdash n$ and $T \in \SYT(\lambda)$, a
  \textit{positive rotation} for $T$ is an interval $[i, k] \subset [n]$
  such that if $T' := (i, i+1, \ldots, k-1, k) \cdot T$, then
  $T' \in \SYT(\lambda)$ and there is some $j$ for which
    \[ \{j\} = \Des(T') - \Des(T)\qquad\text{ and }\qquad
        \{j-1\} = \Des(T) - \Des(T'). \]
\end{Definition}

\noindent
Intuitively, a positive rotation is one for which $j-1 \in \Des(T)$
becomes $j \in \Des(T')$ and all other entries remain the same.
Consequently, $\maj(T') = \maj(T)+1$.  We call $j$ the \emph{moving
  descent} for the positive rotation.

\begin{figure}[ht]
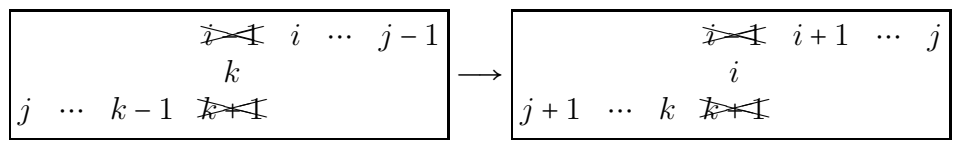
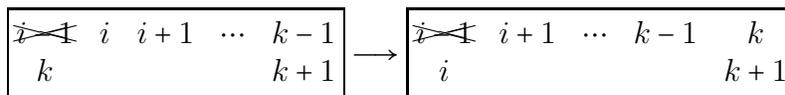
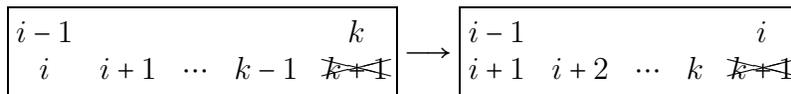

  \centering
  \begin{subfigure}[t]{\textwidth}
    \centering
    \[
      \boxed{\begin{matrix}
        & & & \xcancel{i-1} & i &\cdots & j-1 \\
        & & & k \\
        j &\cdots & k-1 & \xcancel{k+1} \\
      \end{matrix}}
      \longrightarrow
      \boxed{\begin{matrix}
        & & & \xcancel{i-1} & i+1 & \cdots & j \\
        & & & i \\
        j+1 & \cdots & k & \xcancel{k+1}\\
      \end{matrix}}
    \]
    \caption{Schematic of a positive rotation with $i < j < k$.}
    \label{fig:posrot_generic}
  \end{subfigure}
  \begin{subfigure}[t]{\textwidth}
    \centering
    \[
      \boxed{\begin{matrix}
        \xcancel{i-1} & i & i+1 & \cdots & k-1 \\
        k & & & & k+1  \\
      \end{matrix}}
      \longrightarrow
      \boxed{\begin{matrix}
        \xcancel{i-1} &  i+1 & \cdots & k-1 &k\\
        i & & & & k+1  \\
      \end{matrix}}
    \]
    \caption{Schematic of a positive rotation with $i < j=k$.}
    \label{fig:posrot_right}
  \end{subfigure}
  \begin{subfigure}[t]{\textwidth}
    \centering
    \[
      \boxed{
      \begin{matrix}
        i-1 & & & & k \\
        i & i+1 & \cdots & k-1 & \xcancel{k+1} \\
      \end{matrix}}
      \longrightarrow
      \boxed{
      \begin{matrix}
        i-1 & & & & i \\
        i+1 & i+2 & \cdots & k & \xcancel{k+1} \\
      \end{matrix}}
    \]
    \caption{Schematic of a positive rotation with $i=j < k$.}
    \label{fig:posrot_left}
  \end{subfigure}
  \caption{Summary diagrams for positive rotations.}
  \label{fig:posrot}
\end{figure}

The positive rotations can be characterized explicitly
as follows. The proof is omitted  since it follows directly from the
pictures in \Cref{fig:posrot}.

\begin{Lemma}\label{lem:posrot}
  An interval $[i, k]$ is a positive rotation for $T \in \SYT(\lambda)$ if and
  only if $i < k$ and there is some necessarily unique moving descent $j$ with
  $1 \leq i \leq j \leq k \leq n$ such that
  \begin{enumerate}[(a)]
    \item $i, \ldots, j-1$ form a horizontal strip, $j-1, j$ form a vertical
      strip, and $j, j+1, \ldots, k$ form a horizontal strip;
    \item if $i < j$, then $i$ appears strictly northeast of $k$ and
      $i-1$ is not in the rectangle bounding $i$ and $k$;
    \item if $i=j$, then $i-1$ appears in the rectangle bounding $i$ and $k$;
    \item if $j < k$, then $k$ appears strictly northeast of $k-1$ and
      $k+1$ is not in the rectangle bounding $k$ and $k-1$; and
    \item if $j=k$, then $k+1$ appears in the rectangle bounding $k$ and $k-1$.
  \end{enumerate}
\end{Lemma}

\noindent
See \Cref{fig:posrot} for diagrams summarizing these conditions.
\medskip

\begin{figure}[ht]
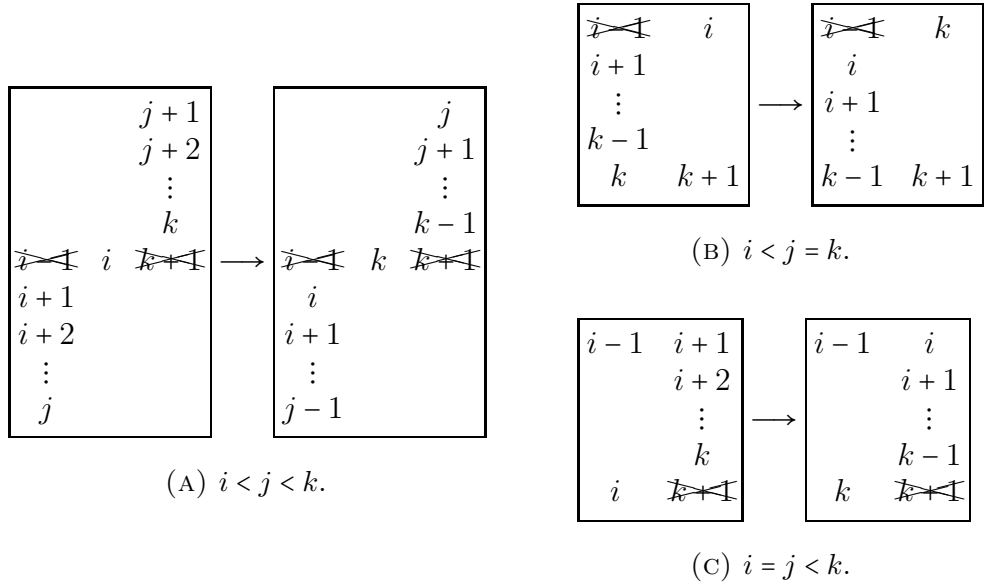

  \centering
  \parbox{\figrasterwd}{
  \parbox{0.55\figrasterwd}{
  \begin{subfigure}[t]{\hsize}
    \centering
    \[
      \boxed{\begin{matrix}
        & & j+1 \\
        & & j+2 \\
        & & \vdots \\
        & & k \\
        \xcancel{i-1} & i & \xcancel{k+1} \\
        i+1 \\
        i+2 \\
        \vdots \\
        j
      \end{matrix}}
      \longrightarrow
      \boxed{\begin{matrix}
        & & j \\
        & & j+1 \\
        & & \vdots \\
        & & k-1 \\
        \xcancel{i-1} & k & \xcancel{k+1} \\
        i \\
        i+1 \\
        \vdots \\
        j-1
      \end{matrix}}
    \]
    \caption{$i < j < k$.}
    \label{fig:negrot_generic}
  \end{subfigure}
  }
  \hskip0.1em
  \parbox{.37\figrasterwd}{
  \begin{subfigure}[t]{\hsize}
    \centering
    \[
      \boxed{\begin{matrix}
        \xcancel{i-1} & i \\
        i+1 \\
        \vdots \\
        k-1 \\
        k & k+1
      \end{matrix}}
      \longrightarrow
      \boxed{\begin{matrix}
        \xcancel{i-1} & k \\
        i \\
        i+1 \\
        \vdots \\
        k-1 & k+1
      \end{matrix}}
    \]
    \caption{$i < j=k$.}
    \label{fig:negrot_right}
  \end{subfigure}
  \vskip1em
  \begin{subfigure}[t]{\hsize}
    \centering
    \[
      \boxed{
      \begin{matrix}
        i-1 & i+1 \\
        & i+2 \\
        & \vdots \\
        & k \\
        i & \xcancel{k+1}
      \end{matrix}}
      \longrightarrow
      \boxed{
      \begin{matrix}
        i-1 & i \\
        & i+1 \\
        & \vdots \\
        & k-1 \\
        k & \xcancel{k+1}
      \end{matrix}}
    \]
    \caption{$i=j < k$.}
    \label{fig:negrot_left}
  \end{subfigure}
  }
  \caption{Summary diagrams for negative rotations.}
  \label{fig:negrot}
  }
\end{figure}

In addition to the \emph{positive rotations} above, we can also
apply \emph{negative rotations}, which are defined exactly as in
\Cref{def:pos_rot} with $(i, i+1, \ldots, k-1, k)$ replaced
by $(k, k-1, \ldots, i+1, i)$ and the rest unchanged.
Combinatorially, negative rotations can be obtained
from positive rotations by applying \emph{inverse-transpose}
moves, that is, by applying negative cycles $(k,k-1,\ldots, i)$
to the transpose of the configurations in \Cref{fig:posrot} and
reversing the arrows. Explicitly, we have the following
analogue of \Cref{lem:posrot}. See \Cref{fig:negrot} for the
cor\-re\-spond\-ing diagrams.

\begin{Lemma}\label{lem:negrot}
  An interval $[i, k]$ is a negative rotation for $T \in \SYT(\lambda)$ if and
  only if $i < k$ and there is some necessarily unique moving descent $j$ with
  $1 \leq i \leq j \leq k \leq n$ such that
  \begin{enumerate}[(a)]
    \item $i, \ldots, j$ form a vertical strip, $j, j+1$ form a horizontal
      strip, and $j+1, \ldots, k$ form a vertical strip;
    \item if $i < j$, then $i+1$ appears strictly southwest of $i$ and
      $i-1$ is not in the rectangle bounding $i$ and $i+1$;
    \item if $i=j$, then $i-1$ appears in the rectangle bounding $i$ and $i+1$;
    \item if $j < k$, then $i$ appears strictly southwest of $k$ and
      $k+1$ is not in the rectangle bounding $i$ and $k$; and
    \item if $j=k$, then $k+1$ appears in the rectangle bounding $i$ and $k$.
  \end{enumerate}
\end{Lemma}

\begin{Example}
  The tableau
  \[
    \begin{matrix}
      1 & 2 & 6 & 7 & 9 \\
      3 & 4 & 8 & 13 & \text{} \\
      5 & 11 & 12 & 15 & \text{} \\
      10 & 14 & \text{} & \text{} & \text{} \\
    \end{matrix}
  \]
  allows positive rotation rules with
  $[i,k] \in \{[5,6], [8,9], [8,10], [8,11],[9,13]\}$, and the tableau
  \[
    \begin{matrix}
      1 & 3 & 8 & 10 & 15 \\
      2 & 4 & 9 & 11 & \text{} \\
      5 & 7 & 13 & 14 & \text{} \\
      6 & 12 & \text{} & \text{} & \text{} \\
    \end{matrix}
  \]
  allows negative rotation rules with
  $[i,k] \in \{[4,6],[6,7],[11,12]\}$.
\end{Example}

\medskip

It turns out that for the vast majority of tableaux, some negative
rotation rule applies.  The positive rotations can be applied in many
of the remaining cases.  For example, among the $81,081$ tableaux in
$\SYT(5442)$, there are only $24$ (i.e., 0.03\%) on which we cannot
apply any positive or negative rotation rule.  For example, no
rotation rules can be applied to the following two tableaux:
\[
  \begin{matrix}
    1 & 2 & 3 & 4 & 5 \\
    6 & 7 & 8 & 9 & \text{} \\
    10 & 11 & 12 & 13 & \text{} \\
    14 & 15 & \text{} & \text{} & \text{} \\
  \end{matrix} \qquad \text{and} \qquad
  \begin{matrix}
    1 & 2 & 3 & 8 & 12 \\
    4 & 6 & 9 & 13 & \text{} \\
    5 & 7 & 10 & 14 & \text{} \\
    11 & 15 & \text{} & \text{} & \text{} \\
  \end{matrix}.
\]

The following lemma and its corollary give a partial explanation for
why negative rotation rules are so common.  Given a tableaux $T$, let $T|_{[z]}$
denote the restriction of $T$ to those values in $[z]$.

\begin{Lemma}\label{lem:negrot.z}
  Let $T \in \SYT(\lambda) \setminus \cE(\lambda)$. Suppose
  $z$ is the largest value such that $T|_{[z]}$ is
  contained in $\maxmaj(\mu)$ for some $\mu$.
  If $T|_{[z+1]}$ is not of the form
    \[ \begin{matrix}
          1 & 2 & \cdots & i \\
          i+1 & z+1 \\
          i+2 &  \\
          \vdots \\
          z
        \end{matrix} \]
  then some negative rotation rule applies to $T$.

  \begin{proof}
    Since $T \not\in \cE(\lambda)$, $T$ is not
    $\maxmaj(\lambda)$, so $\lambda$ is not
    a one row or column shape. We have $z \geq 2$
    since both two-cell tableaux are the max-maj tableau of
    their shape. Since $\maxmaj(\mu)$ is built from
    successive, outermost, maximal length, vertical strips as in \Cref{fig:maxmaj_ex},
    the same is true of $T|_{[z]}$.

    First, suppose $z$ is not in the lowest row of $T|_{[z]}$.
    Let $i$ be the value in the topmost corner cell in $T|_{[z]}$ which is
    strictly below $z$. Let $j \geq i$ be the bottommost
    cell in the vertical strip of $T|_{[z]}$ which contains $i$.
    See \Cref{fig.z.not.low}.
    We verify the conditions of \Cref{lem:negrot}, so the
    negative $[i, z]$-rotation rule applies with moving descent $j$.
    By construction, $i, \ldots, j$ form a vertical strip,
    $j, j+1$ form a horizontal strip, and $j+1, \ldots, z$
    form a vertical strip. If $i<j$, then since $i$ is a corner
    cell, $i+1$ appears strictly southwest of $i$, and
    $i-1$ is above both $i$ and $i+1$ so $i-1$ is not in the
    rectangle bounding $i$ and $i+1$. If $i=j$, we see that
    $i-1$ appears in the rectangle bounded by $i$ and $i+1$.
    We also see that $i$ appears strictly southwest of $z$,
    and $z+1$ is not in the rectangle bounding $i$ and $z$
    since $i$ is a topmost corner and $z$ is maximal.

\begin{figure}[ht]
  \centering
  \begin{subfigure}[t]{\textwidth}
    \centering
      $$
              \begin{matrix}
          1&3&6&\red{11}\\
          2&4&7&\red{12}\\
          5&\red{8}&&\\
          \red{9}&13&&\\
          \red{10}&&&\\
        \end{matrix}
        \hspace{.2in}        \longrightarrow \hspace{.2in}
                      \begin{matrix}
          1&3&6&\red{10}\\
          2&4&7&\red{11}\\
          5&\red{12}&&\\
          \red{8}&13&&\\
          \red{9}&&&\\
        \end{matrix}
        $$
        \caption{For the tableau on the left above,  $i=8$ and $z=12$ since $T|_{[12]}$ is contained the
          max-maj tableau of shape 44322, 12 is not in the lowest
          row, $8$ is in the closet corner to 12 in $T|_{[12]}$ and
          below 12.  Apply the negative rotation $(12,11,10,9,8)$
                    to get the tableau on the right, and observe maj has
          increased by 1.  The moving descent is $j=10$.}\label{fig.z.not.low}
  \end{subfigure}
  \begin{subfigure}[t]{\textwidth}
    \centering
      $$
              \begin{matrix}
          1&3&6\\
          2&4&\red{7}\\
          5&\red{8}&11\\
          \red{9}&&\\
          \red{10}&&\\
        \end{matrix}
\hspace{.2in}        \longrightarrow \hspace{.2in}
              \begin{matrix}
          1&3&6\\
          2&4&\red{10}\\
          5&\red{7}&11\\
          \red{8}&&\\
          \red{9}&&\\
        \end{matrix}
        $$
        \caption{For the tableau on the left above,  $i=7$ and $z=10$ since $T|_{[10]}$ is the
          max-maj tableau of shape 33211, 10 is in the lowest
          row, $11$ is in row 3, and $7$ is the largest value in $T|_{[10]}$
          in row 2.  Apply the negative rotation $(10,9,8,7)$
          to get the tableau on the right, and observe maj has
          increased by 1.  The moving descent is $j=z=10$.}
        \label{fig.z.low}
  \end{subfigure}
  \caption{Examples of the negative rotations obtained from
    \Cref{lem:negrot.z}.}
  \label{fig:z.exs}
\end{figure}

    Now suppose $z$ is in the lowest row of $T|_{[z]}$.  In this case,
    $T|_{[z]}$ is the max-maj tableau of its shape, so that
    $z < |\lambda|$ and $z+1$ exists in $T$ since
    $T \not\in \cE(\lambda)$.  By maximality of $z$, $z+1$ cannot be
    in row $1$ or below $z$. Let $i<z$ be the value in the
    rightmost cell of $T|_{[z]}$ in the row immediately above
    $z+1$. See \Cref{fig.z.low}.  We check that the negative
    $[i, z]$-rotation rule applies with moving descent $j=z$ using the
    conditions in \Cref{lem:negrot}. By construction, $i, \ldots, z$
    form a vertical strip. Since $z+1$ is not below $z$, we see that
    $z, z+1$ form a horizontal strip.  Since $z+1$ is in the row
    below $i$, $i+1$ appears strictly southwest of $i$.  We also see
    that $z+1$ appears in the rectangle bounded by $i$ and $z$ by
    choice of $i$.  It remains to show that $i-1$ is not in the
    rectangle bounding $i$ and $i+1$.  Suppose to the contrary that
    $i-1$ is in the rectangle bounding $i$ and $i+1$.  Then $i$ would
    have to be in row $1$ by choice of $i<z$.  Consequently $i+1$ is
    in row $2$ and strictly west of $i$, forcing $i-1$ to be in row
    $1$ also. It follows from the choice of $z$ that $T|_{[i]}$ is a
    single row, the values $i,i+1,\ldots z$ form a vertical strip, and
    $T|_{[z+1]}$ is of the above forbidden form, giving a
    contradiction.
  \end{proof}
\end{Lemma}

\begin{Corollary}\label{cor:descent.at.1}
  If $T \in \SYT(\lambda) \setminus \cE(\lambda)$ and
  $1 \in \Des(T)$, then some negative rotation rule
  applies to $T$.

  \begin{proof}
    Let $z$ be as in \Cref{lem:negrot.z}. Clearly
    $z \geq 2$ and $T|_{[2]}$ is a single column,
    so $T|_{[z+1]}$ cannot possibly be of the forbidden form.
  \end{proof}
\end{Corollary}


We also have the following variation on \Cref{lem:negrot.z}.  It is
based on finding the largest value $q$ such that $T|_{[q]}$ is
contained in an exceptional tableau of type (iii).  The proof is again
a straightforward verification of the conditions in \Cref{lem:negrot},
and is omitted.

\begin{Lemma}\label{lem:negrot.rect.except}
  Let $T \in \SYT(\lambda) \setminus \cE(\lambda)$. Suppose
the initial values of $T$ are of the form
$$
    \begin{matrix}
      1 & 2\\
      3 & p+1\\
      4 & \vdots \\
      \vdots & q\\
      \vdots &  \xcancel{q+1}\\
      p &  \\
    \end{matrix}
    \hspace{1cm}
    \text{or}
    \hspace{1cm}
        \begin{matrix}
      1 & 2 & \ell+1 & \cdots&  \vdots & p+1 \\
      3 & z+1 & \vdots & \vdots & \vdots &\vdots \\
      4 & z+2 & \vdots & \vdots & \vdots & q\\
      \vdots & \vdots  & \vdots & \vdots &  \vdots & \xcancel{q+1} \\
      z & \ell & m &\cdots& p & \\
    \end{matrix}.
    $$
    In either case, the $[p,q]$-negative rotation rule
    applies to $T$.
    \end{Lemma}
\subsection{Initial Block Rules}

Here we describe a collection of five additional \textit{block rules} which may apply
to a tableau that is not in the exceptional set.  In each case, if the
rule applies, then we specify a permutation of the entries so that we
either add $1$ into the descent set and leave the other descents
unchanged, or we add $1$ into the descent set, increase one existing descent by
$1$, and decrease one existing descent by $1$.  Thus, maj will increase by
1 in all cases. While these additional rules are certainly not uniquely determined
by these criteria, they are also not arbitrary.

\begin{Example}
  For a given $T \in \SYT(\lambda)$, one may consider all
  $T' \in \SYT(\lambda)$ where $\maj(T') = \maj(T) + 1$.
  If $T' = \sigma \cdot T$ where $\sigma$ is a simple cycle,
  then one of the rotation rules may apply to $T$.
  \Cref{tab:B_examples} summarizes five particular $T$ for
  which \textit{no} rotation rules apply. These examples have
  guided our choices in defining the block rules. In all but one
  of these examples, there is a unique $T'$ with
  $\maj(T') = \maj(T)+1$, though in the
  third case there are two such $T'$, one of which ends up
  being easier to generalize.
  \begin{table}[ht]
  \begin{center}
  \begin{tabu}{c|c|c|c}
    Tableau $T$
      & Tableaux $T'$
      & $\sigma$
      & Block rule \\
      \hline\hline
    $\begin{smallmatrix}
      1 & 2 & 3 & 7 \\
      4 & 5 & 6 & 8 \\
    \end{smallmatrix}$
      & $\begin{smallmatrix}
          1 & 3 & 4 & 6 \\
          2 & 5 & 7 & 8 \\
          \end{smallmatrix} $
      & $(2,3,4)(6,7)$
      & B1 \\
    \hline
    $\begin{smallmatrix}
      1 & 2 & 3 & 4 \\
      5 & 6 & 7 \\
    \end{smallmatrix}$
      & $\begin{smallmatrix}
          1 & 3 & 4 & 7 \\
          2 & 5 & 6 &  \\
          \end{smallmatrix}$
      & $(2,3,4,7,6,5)$
      & B2 \\
    \hline
    $\begin{smallmatrix}
      1 & 2 & 3 \\
      4 & 6 &  \\
      5 & 7 &  \\
    \end{smallmatrix}$
      & $\begin{smallmatrix}
          1 & 3 & 6 \\
          2 & 4 &  \\
          5 & 7 &  \\
          \end{smallmatrix},
          \begin{smallmatrix}
          1 & 4 & 5 \\
          2 & 6 &  \\
          3 & 7 &  \\
          \end{smallmatrix}$
      & $(2,3,6,4), (2,4)(3,5)$
      & B3, --- \\
    \hline
    $\begin{smallmatrix}
      1 & 2 & 7 \\
      3 & 5 & 8 \\
      4 & 6 & 9 \\
      10 & &
    \end{smallmatrix}$
      & $\begin{smallmatrix}
          1 & 4 & 8 \\
          2 & 5 & 9 \\
          3 & 6 & 10 \\
          7 & &
          \end{smallmatrix}$
      & $(2,4,3)(7,8,9,10)$
      & B4 \\
    \hline
    $\begin{smallmatrix}
    1 & 2 \\
    3 & 5 \\
    4 & 6 \\
    7 &
    \end{smallmatrix}$
      & $\begin{smallmatrix}
          1 & 5 \\
          2 & 6 \\
          3 & 7 \\
          4 &
          \end{smallmatrix}$
      & $(2,5,6,7,4,3)$
      & B5
  \end{tabu}
  \end{center}
  \caption{Some tableaux $T \in \SYT(\lambda)$ together with all
    $T' = \sigma \cdot T \in \SYT(\lambda)$ where
    $\maj(T') = \maj(T)+1$. See \Cref{def:block_rules} for an
    explanation of the final column.}\label{tab:B_examples}
  \end{table}
\end{Example}

\medskip

In the remainder of this subsection, we describe the block rules,
abbreviated B-rules.  Then, we prove that if no rotation rules are
possible for a tableau then either it is in the exceptional set or we
can apply one of the B-rules.  The B-rules cover disjoint cases so no
tableau admits more than one block rule.  To state the B-rules
precisely, assume $T\in \SYT(\lambda) \setminus \cE(\lambda)$
and no rotation rule applies.

\begin{Notation}\label{not:abc}
  Let $c$ be the largest possible value such that $T|_{[c]}$ is
  contained in the min-maj tableau of a rectangle shape with $a$
  columns and $b$ rows.  Consequently, the first $a$ numbers in row
  $i$, $1 \leq i \leq b-1$, of $T$ are $(i-1)a + 1,\ldots, ia$, and
  row $b$ begins with $(b-1)a+1, (b-1)a+2, \ldots, c$.
\end{Notation}

Assuming $1 \not \in \Des(T)$ and $T \not \in\cE(\lambda)$, we know
$a,b\geq 2$ and $c \geq 3$.  If $c+1$ is in $T$, then it must be
either in position $(1,a+1)$ or $(b+1,1)$.  If $c=ab$, then
$c < |\lambda|$ since $T \not \in\cE(\lambda)$, otherwise
$c = |\lambda|$ is possible.  For example, the tableaux
$$\begin{smallmatrix}
1 & 2 & 3 & 4 & 5 & 16 \\
6 & 7 & 8 & 9 &10 & 17 \\
11 & 12 & 13 & 14 & 15 &
\end{smallmatrix}, \quad \begin{smallmatrix}
1 & 2 & 3 & 4 & 5 \\
6 & 7 & 8 & 9 & 10 \\
11 & 12 & 13 & &
\end{smallmatrix}, \quad \begin{smallmatrix}
1 & 2 & 3 & 4 \\
5 & 9 & & \\
6 & 10 & & \\
7 & & & \\
8 & & &
\end{smallmatrix}, \quad \begin{smallmatrix}
1 & 2 & 7 & 10 \\
3 & 5 & 8 & 11 \\
4 & 6 & 9 & 12 \\
13 & & &
\end{smallmatrix}
$$
have $(a,b,c)$ equal to $(5,3,15)$, $(5,3,13)$, $(4,2,5)$, and $(2,2,3)$, respectively.

\medskip

\begin{Definition}\label{def:block_rules}
  Using \Cref{not:abc}, we identify the \emph{block
    rules} with further required assumptions as follows. See
    \Cref{fig:Brules} for summary diagrams.
\end{Definition}
\begin{itemize}
\item \textbf{Rule B1}: Assume $c=ab$, $T_{(1,a+1)}=c+1$,
  $T_{(2,a+1)}= c+2$, and $a<c-2$.  In this case, we perform the
  rotations $(2,\ldots,a+1)$ and $(c,c+1)$ which are sufficiently
  separated by hypothesis.  Then, $1,a+1$ and $c$ become descents, and
  $a$ and $c+1$ are no longer descents, so the major index is
  increased by $1$. The B1 rule is illustrated here with $a = 5$,
  $b = 3$:

$$\mbox{B1:} \qquad  \begin{matrix}
  1 & \red 2 & \red 3 & \red 4 & \boxed{\red 5} & \boxed{\blue{16}}\\
  \red 6 & 7 & 8 & 9 &10 & 17\\
  11 & 12 & 13 & 14 & \blue{15} &
  \end{matrix}  \quad \begin{matrix}\red \ccwise \\ \blue \ccwise \end{matrix} \quad
  \begin{matrix}
    \boxed 1 & \red 3 & \red 4 & \red 5 & \boxed{\red 6} & \boxed{\blue{15}}\\
    \red 2 & 7 & 8 & 9 &10 & 17\\
    11 & 12 & 13 & 14 & \blue{16} &
  \end{matrix}$$

\medskip

The boxed numbers represent descents of the tableau on the left/right that are not descents of the tableau on the right/left. The elements not shown (i.e. $18,19,\ldots,|\lambda|$) can be in any position.

\medskip

\item \textbf{Rule B2}: Assume $c<ab$ and there exists a $1 \leq k<a$
  such that $T_{(b,k)}=c$ and $T_{(b,k+1)}\neq c+1$.  In this case, we
  perform the rotation
  $(2,3,\ldots,a,2a,3a,\ldots,a(b-1),c,c-1,\ldots,c-k+1=a(b-1)+1,a(b-2)+1,\ldots,2a+1,a+1)$
  around the perimeter of $T|_{[c]}$. Now $1$ becomes a descent, and
  the other descents stay the same so the major index again increases
  by $1$. The B2 rule is illustrated by the following (here $a = 5$,
  $b = 2$ and $k = 3$):

$$  \mbox{B2:} \qquad  \begin{matrix}
  1 & \red 2 & \red 3 & \red 4 & \red 5 \\
  \red 6 & 7 & 8 & 9 & \red{10} \\
  \red{11} & \red{12} & \red{13} & \xcancel{14}&
  \end{matrix}  \quad \begin{matrix}\red \ccwise \end{matrix} \quad
  \begin{matrix}
  \boxed 1 & \red 3 & \red 4 & \red 5 & \red{10} \\
  \red 2 & 7 & 8 & 9 & \red{13}\\
  \red 6 & \red{11} & \red{12} & \xcancel{14}&
\end{matrix}$$
The crossed out number $14$ means that $14$ is not in position $(3,4)$: it can either be in positions $(1,6)$ or $(4,1)$, or it can be that $\lambda = 553$. Again, the numbers $15,\ldots,|\lambda|$ can be anywhere in $T$.

\medskip

\item
\textbf{Rule B3}: Assume $a\geq 3$, $c=a+1$, and there exists
$k \geq 2$ such that $T_{(2,2)} = a+k+1$, $T_{(3,2)} = a+k+2$, and for all $i
\in \{ 1,2,\ldots,k\}$ we have
$T_{(i+1,1)}=a+i$.  Thus $b=2$. Then we apply the rotation
$(2,3,\ldots,a,a+k+1,a+1)$. Now
$1$ becomes a descent, and the rest of the descent set is unchanged so
the major index again increases by
$1$. The B3 rule is illustrated by the following (here $a = 4$, $k =
4$):

$$\mbox{B3:} \qquad  \begin{matrix}
  1 & \red 2 & \red 3 & \red 4  \\
  \red 5 & \red 9 & & \\
  6 & 10 & & \\
  7 & & &\\
  8 & & &
  \end{matrix}  \quad \begin{matrix}\red \ccwise \end{matrix}  \quad
  \begin{matrix}
  \boxed 1 & \red 3 &  \red 4 & \red 9 \\
  \red 2 & \red 5 & & \\
  6 & 10 & & \\
  7 & & &\\
  8 & & &
  \end{matrix}$$

\medskip

\item \textbf{Rule B4}: Assume that $a = 2$, $c=3$, and there exists
  $k \geq 2$ such that $\{3,4,\ldots,k+1\}$ appear in column 1 of $T$,
  $\{k+2, k+3,\ldots, 2k\}$ appear in column 2 in $T$.  Further assume
  that the set $\{2k+1,2k+2,\ldots,
 3k\}$ appears in column 3, $\{3k+1,3k+2,\ldots, 4k\}$ appears in column 4, etc., until
  column $l$ for some $l>2$ and $T_{(k+1,1)} = kl+1$ and
  $T_{(k+1,2)} \neq kl+2$.  Thus, $b=2$. In this case, we can perform the two
  rotations $(k+1,k,\ldots,3,2)$ and
  $(k(l-1)+1,k(l-1)+2,\ldots,kl,kl+1)$. Now $1$, $k+1$ and $k(l-1)$
  enter the descent set, and $k$ and $k(l-1)+1$ leave it, so the major
  index increases by $1$.  The B4 rule is illustrated by the following
  (here $k = 3$ and $l = 4$):
  $$
  \mbox{B4:} \qquad  \begin{matrix}
  1 & {\red 2} & 7 & \boxed{\blue {10}} \\
  \boxed{\red 3} & 5 & 8 &\blue{11}\\
  \red 4 & 6 & 9 &\blue{12}\\
  \blue{13} & \xcancel{14} & &
  \end{matrix}  \quad\begin{matrix}\red \ccwise \\ \blue \ccwise \end{matrix}  \quad
  \begin{matrix}
  \boxed 1 & \boxed{\red 4} & 7 & \blue{11} \\
  \red 2 & 5 & 8 &\blue{12}\\
  \red 3 & 6 & \boxed 9 &\blue{13}\\
  \blue{10} & \xcancel{14} & &
  \end{matrix}$$

\item
\textbf{Rule B5}: Assume that $a = 2$, $c=3$, and there exists
$k > 3$ such that $\{3,4,\ldots,k\}$ appear in column 1 of $T$,
$\{k+1, k+2,\ldots, 2k-2\}$ appear in column 2 in $T$.  Furthermore,
assume $T_{(k,1)} = 2k-1$ and $T_{(k,2)} \neq 2k$. Thus, $b=2$. Then apply the
cycle $(k,k-1,\ldots,3,2,k+1,k+2,\ldots,2k-1)$ to $T$. Now $1$ becomes
a descent, and the rest of the descent set remains unchanged, so the
major index increases by $1$.  The B5 rule is illustrated by the
following (here $k = 5$):
  $$
  \mbox{B5:} \qquad \begin{matrix}
  1 & {\red 2}  \\
  {\red 3} & \red 6\\
  \red 4 & \red 7\\
  \red 5 & \red 8 \\
  \red 9 & \xcancel{10}
  \end{matrix}  \quad \begin{matrix}\red \ccwise \end{matrix} \quad
  \begin{matrix}
  \boxed 1 & \red 6  \\
  \red 2 & \red 7\\
  \red 3 & \red 8\\
  \red 4 & \red 9 \\
  \red 5 & \xcancel{10}
\end{matrix}$$

\end{itemize}

\begin{figure}[b]
  \centering
  \begin{subfigure}[t]{\textwidth}
    \centering
    \begin{gather*}
      \boxed{\begin{matrix}
        1 & \red{2} & \cdots & \red{a} & \blue{ab+1} \\
        \red{a+1} & a+2 & \cdots & 2a & ab+2 \\
        \vdots & \vdots & \ddots & \vdots \\
        a(b-1)+1 & a(b-1)+2 & \cdots & \blue{ab}
      \end{matrix}}\\
      \downarrow\\
      \boxed{\begin{matrix}
        1 & \red{3} & \cdots & \red{a+1} & \blue{ab} \\
        \red{2} & a+2 & \cdots & 2a & ab+2 \\
        \vdots & \vdots & \ddots & \vdots \\
        a(b-1)+1 & a(b-1)+2 & \cdots & \blue{ab+1}
      \end{matrix}}
    \end{gather*}
    \caption{B1.}
    \label{fig:Brules.1}
  \end{subfigure}
  \begin{subfigure}[t]{\textwidth}
    \centering
    \begin{gather*}
      \boxed{\begin{matrix}
        1 & \red{2} & \cdots & \cdots & \cdots & \cdots & \red{a} \\
        \red{a+1} & a+2 & \cdots & \cdots & \cdots & \cdots & \red{2a} \\
        \vdots & \vdots & \ddots & \vdots & \vdots & \ddots & \vdots \\
        \red{a(b-2)+1} & a(b-2)+2 & \cdots & \cdots & \cdots & \cdots & \red{a(b-1)} \\
        \red{a(b-1)+a} & \red{a(b-1)+2} & \cdots & \red{c-1} & \red{c} & \cdots & \xcancel{ab}
      \end{matrix}}\\
      \downarrow\\
      \boxed{\begin{matrix}
        1 & \red{3} & \cdots & \cdots & \cdots & \cdots & \red{2a} \\
        \red{2} & a+2 & \cdots & \cdots & \cdots & \cdots & \red{3a} \\
        \vdots & \vdots & \ddots & \vdots & \vdots & \ddots & \vdots \\
        \red{a(b-3)+1} & a(b-2)+2 & \cdots & \cdots & \cdots & \cdots & \red{c} \\
        \red{a(b-2)+1} & \red{a(b-1)+1} & \cdots & \red{c-2} & \red{c-1} & \cdots & \xcancel{ab}
      \end{matrix}}
    \end{gather*}
    \caption{B2.}
    \label{fig:Brules.2}
  \end{subfigure}
  \begin{subfigure}[t]{\textwidth}
    \centering
    \[
      \boxed{\begin{matrix}
        1 & \red{2} & \cdots & \red{a-1} & \red{a} \\
        \red{a+1} & \red{a+k+1} \\
        a+2 & a+k+2 \\
        \vdots \\
        a+k
      \end{matrix}}
      \longrightarrow
      \boxed{\begin{matrix}
        1 & \red{3} & \cdots & \red{a} & \red{a+k+1} \\
        \red{2} & \red{a+1} \\
        a+2 \\
        \vdots \\
        a+k
      \end{matrix}}
    \]
    \caption{B3.}
    \label{fig:Brules.3}
  \end{subfigure}
\end{figure}
\begin{figure}[t]\ContinuedFloat
  \begin{subfigure}[t]{\textwidth}
    \centering
    \begin{gather*}
      \boxed{\begin{matrix}
        1 & \red{2} & 2k+1 & \cdots & \blue{k(\ell-1)+1} \\
        \red{3} & k+2 & 2k+2 & \cdots & \blue{k(\ell-1)+2} \\
        \red{4} & k+3 & 2k+3 & \cdots & \blue{k(\ell-1)+3} \\
        \vdots & \vdots & \vdots & \ddots & \vdots \\
        \red{k+1} & 2k & 3k & \cdots & \blue{k\ell} \\
        \blue{k\ell+1} & \xcancel{k\ell+2}
      \end{matrix}}\\
      \downarrow\\
      \boxed{\begin{matrix}
        1 & \red{k+1} & 2k+1 & \cdots & \blue{k(\ell-1)+2} \\
        \red{2} & k+2 & 2k+2 & \cdots & \blue{k(\ell-1)+3} \\
        \red{3} & k+3 & 2k+3 & \cdots & \blue{k(\ell-1)+4} \\
        \vdots & \vdots & \vdots & \ddots & \vdots \\
        \red{k} & 2k & 3k & \cdots & \blue{k\ell+1} \\
        \blue{k(\ell-1)+1} & \xcancel{k\ell+2}
      \end{matrix}}
    \end{gather*}
    \caption{B4.}
    \label{fig:Brules.4}
  \end{subfigure}
  \begin{subfigure}[t]{\textwidth}
    \centering
    \[
      \boxed{\begin{matrix}
        1 & \red{2} \\
        \red{3} & \red{k+1} \\
        \red{4} & \red{k+2} \\
        \vdots & \vdots \\
        \red{k-1} & \red{2k-3} \\
        \red{k} & \red{2k-2} \\
        \red{2k-1} & \xcancel{2k}
      \end{matrix}}
      \longrightarrow
      \boxed{\begin{matrix}
        1 & \red{k+1} \\
        \red{2} & \red{k+2} \\
        \red{3} & \red{k+3} \\
        \vdots & \vdots \\
        \red{k-2} & \red{2k-2} \\
        \red{k-1} & \red{2k-1} \\
        \red{k} & \xcancel{2k}
      \end{matrix}}
    \]
    \caption{B5.}
    \label{fig:Brules.5}
  \end{subfigure}
  \caption{Summary diagrams for block rules.}
  \label{fig:Brules}
\end{figure}

\begin{Lemma}\label{lem:descent.after.2}
  If $T\in \SYT(\lambda)$, \ $T\not \in \cE(\lambda)$, and
  $1,2 \not \in \Des(T)$, then either some rotation rule applies to $T$
  or a B1, B2 or B3 rule applies.
\end{Lemma}

\begin{proof}
  Let $c$ be the largest possible value such that $T|_{[c]}$ is
  contained in the min-maj tableau of a rectangle shape with $a$
  columns and $b$ rows, see \Cref{not:abc}.
  Since $1,2 \not \in \Des(T)$ and $T\not \in \cE(\lambda)$, we
  know $1,2,3$ are in the first row  of $T$ so $a\geq 3$, $b \geq 2$,
  and $a+2\leq |\lambda|$.  By construction, we have $T_{(2,1)}=a+1$
  and $a+2$ must appear in position $(1,a+1)$, $(2,2)$, or $(3,1)$ in
  $T$.

 \medskip

 \noindent
 \textbf{Case 1: $T_{(1,a+1)}=a+2$.}  Observe that
  $$
T|_{[a+2]}\ =  \   \begin{matrix}
      1 & 2 &  3 & \cdots & a & a+2  \\
      a+1 &  & &  & \\
  \end{matrix}
  $$
  and $z \geq a+2$. Consequently, $T|_{[z+1]}$ cannot
  be of the form forbidden by \Cref{lem:negrot.z}, so
  a negative rotation rule applies.

  \medskip
  \noindent
  \textbf{Case 2: $T_{(2,2)}=a+2$.}
  First suppose $c=ab$, then $T_{(1,a+1)}=c+1$ by choice of $c$.  Now consider
  the two subcases, $T_{(2,a+1)}=c+2$ and $T_{(2,a+1)}\neq c+2$.  In
  the former case, as in \Cref{fig:Brules.1}, the B1 rule applies to $T$. In the
  latter case, one may check that an $[i,c+1]$-positive rotation rule applies to $T$ where
  $i=T_{(b,1)}$. On the other hand, if $c< ab$, then a B2 rule applies to $T$
  as in \Cref{fig:Brules.2}.

%

\noindent
\textbf{Case 3: $T_{(3,1)}=a+2$.}  Let
$k= \min \{j \geq 2 \ | \ a+j \not \in \Des(T)\}$ so $T_{(k+1,1)}=a+k$
and $T_{(k+2,1)}\neq a+k+1$.  Since $T \not \in \cE(\lambda)$,
we know $a+k+1$ exists in $T$ either in position $(1,a+1)$ or $(2,2)$,
so $T|_{[a+k+1]}$ looks like
$$
   \begin{matrix}
      1 & 2 &  3 & \cdots & a  & a+k+1\\
      a+1 &  & &  \\
      a+2 &  & &  \\
      \vdots &  & &  \\
      a+k &  & &  \\
  \end{matrix}
\hspace{.4in} \text{  or  } \hspace{.4in}
  \begin{matrix}
      1 & 2 &  3 & \cdots & a \\
      a+1 & a+k+1 & &  \\
      a+2 &  & &  \\
      \vdots &  & &  \\
      a+k &  & &  \\
    \end{matrix} \, .
$$
If $T_{(1,a+1)}=a+k+1$, then \Cref{lem:negrot.z} shows that a negative
rotation rule applies to $T$. On the other hand, if
$T_{(2,2)}=a+k+1$, then observe that either a B3 move applies or the
rotation $(a+k,a+k+1)$ applies to $T$, depending on whether
$T_{(3,2)}=a+k+2$ or not.
\end{proof}

\begin{Lemma}\label{lem:descent.at.2}
  If $T\in \SYT(\lambda)$, \ $T\not \in \cE(\lambda)$,
  $1\not \in \Des(T)$, and $2 \in \Des(T)$, then either some rotation rule
  applies to $T$ or a B1, B2, B4 or B5 rule applies.
\end{Lemma}

\begin{proof}
  Let $k = \min \{j \geq 3\ | \ j \not \in \Des(T)\}$ so the consecutive
  sequence $[3,k]$ appears in the first column of $T$ and $k+1$ does
  not.  By definition of $k$ and the fact that
  $T\not \in \cE(\lambda)$, $T$ must have $k+1$ in position
  $(1,3)$ or $(2,2)$.  If $T_{(1,3)}=k+1$, then a negative rotation
  rule holds by \Cref{lem:negrot.z}.

  Assume $T_{(2,2)}=k+1$.  Let $\ell$ be the maximum value such that
  $[k+1,\ell]$ appears as a consecutive sequence in column 2 of $T$.
  If $\ell < 2(k-1)$, then the negative rotation rule for
  $(\ell, \ell-1,\ldots, k)$ applies to $T$ by the first
  case of \Cref{lem:negrot.rect.except}.

  If $\ell = 2(k-1)$ and $T_{(1,3)}=\ell+1$, let
  $m$ be the maximum value such that
  $[\ell+1,m]$ appears as a consecutive sequence in column 3 of
  $T$.  We subdivide on cases for $m$ again.  If
  $m<3(k-1)$, then the negative rotation rule
  $(m,m-1,\ldots,\ell)$ applies to
  $T$ by the second case of \Cref{lem:negrot.rect.except}. If
  $m=3(k-1)$, we consider the maximal sequence of columns containing a
  consecutive sequence in rows
  $[1,k-1]$ to the right of column 2 until one of two conditions hold
  $$
    \boxed{
    \begin{matrix}
      1 & 2 & \ell+1 & \cdots& \vdots \\
      3 & k+1 & \vdots & \vdots & p\\
      \vdots & \vdots  & \vdots & \ddots &  \xcancel{p+1} \\
      k & \ell & m &\cdots& \\
    \end{matrix}}
  \hspace{1in}
      \boxed{
    \begin{matrix}
      1 & 2 & \ell+1 & \cdots& \\
      3 & k+1 & \vdots & \vdots & \vdots\\
      \vdots & \vdots  & \vdots & \ddots & \vdots \\
      k & \ell & m & \cdots & p\\
      p+1 &&& \\
    \end{matrix}}
  $$
  In the first picture, $T|_{[p]}$ is not a rectangle, so we may apply a
  negative rotation by the second case of \Cref{lem:negrot.rect.except},
  so consider the second picture.
  In the second picture, $T|_{[p]}$ is a rectangle and we know $p+1$
  exists in $T$ since $T|_{[p]}$ is an exceptional tableau for a
  rectangle shape. If $p+2$ is in row $k$, column $2$, a
  $(p,p+1)$ rotation rule applies.  If $p+2$ is not in row $k$, column $2$,
  then a B4-move applies.

  Finally, consider the case $\ell = 2(k-1)$ and $T_{(k,1)}=\ell+1$.
 If $T_{(k,2)} \neq \ell+2$ and $k>3$, then a B5-move applies.  If
$T_{(k,2)}=\ell+2$ and $k>3$, then the rotation $(\ell,\ell+1)$ applies to $T$
since $\ell-1$ is above $\ell$.
  If $T_{(k,2)}=\ell+2$ and $k=3$, then $\ell=4=T_{(2,2)}$ and
    $T_{(3,1)}=5$ so  $T$ contains
  $$
    \begin{matrix}
      1 & 2 \\
      3 & 4 \\
      5 &
    \end{matrix} \ .
  $$
  In this case, consider the subcases  $c=ab$ or $c<ab$ with $a=2$.  If $c=ab$, then
  $T_{(1,3)}=c+1$ since $T\not \in \cE(\lambda)$.  Either a
  B1-move applies if $T_{(2,3)}=c+2$ and a $(c,c+1)$ rotation applies
  otherwise.  On the other hand, if $c<ab$ then a B2-rule applies.
\end{proof}

We may finally define the map $\varphi$ from
\eqref{eq:varphi_telegraph}.  The proof of \Cref{thm:zeros} from the
introduction follows immediately from this definition and the last few
lemmas.

\begin{Definition}
  Given $T \in \SYT(\lambda) - \cE(\lambda)$, we define $\varphi(T)$
  as follows. If $1 \in \Des(T)$, define
  $\varphi(T) = (z,z-1,\ldots,i)T$ as in \Cref{cor:descent.at.1}.  If $1,2 \not \in \Des(T)$, then
  \Cref{lem:descent.after.2} applies, so define $\varphi(T)$
  using the specific B1, B2, B3 or rotation rule identified in the
  proof of that lemma.  If $1 \not \in \Des(T)$ and $2 \in \Des(T)$, then
  \Cref{lem:descent.at.2} applies, so define $\varphi(T)$ using the
  specific B1, B2, B4, B5, or negative rotation rule identified in the proof
  of that lemma. These rules cover all possible cases. By contruction,
  $\maj(\varphi(T)) = \maj(T)+1$.
\end{Definition}

\medskip

We may define two poset structures on standard tableaux of a given
shape using the preceding combinatorial operations. We call
them ``strong'' and ``weak'' in analogy with the strong and weak
Bruhat orders on permutations. Recall an \textit{inverse-transpose}
block rule is a block rule obtained from transposing the diagrams in
\Cref{fig:Brules} and reversing the arrows.

\begin{Definition}
  As sets, let $P(\lambda)$ and $Q(\lambda)$ be either
    \[ \SYT(\lambda) \setminus \{\minmaj(\lambda), \maxmaj(\lambda) \} \]
  if $\lambda$ is a rectangle with at least two rows and columns, or
  $\SYT(\lambda)$ otherwise.
  \begin{itemize}
    \item (Strong SYT Poset) Let $P(\lambda)$ be the partial
      order with covering re\-la\-tions given by rotations, block rules, and
      inverse-transpose block rules increasing $\maj$ by $1$.
    \item (Weak SYT Poset) Let $Q(\lambda)$ be the
      partial order with covering relations given by
      $S \prec T$ if $\varphi(S)=T$ or
      $\varphi(T')=S'$ where $S',T'$ are the transpose of
      $S,T$, respectively.
  \end{itemize}
\end{Definition}

\begin{Corollary}
  As posets, $P(\lambda)$ and $Q(\lambda)$ are ranked
  with a unique minimal and maximal element. If $\lambda$ is
  not a rectangle, the rank function is given by $\rk(T)=\maj(T)-b(\lambda)$.
  If $\lambda$ is a rectangle with at least two rows and columns, then the
  rank function is given by $\rk(T)=\maj(T)-b(\lambda)-2$.

  \begin{proof}
    By \Cref{cor:first_few_coeffs}, $P(\lambda)$ and
    $Q(\lambda)$ have a single element of minimal $\maj$
    and of maximal $\maj$. Any element $T$ besides these
    is covered by $\varphi(T)$ and covers $\varphi(T')'$, so
    is not maximal or minimal. By construction $\maj$ increases
    by $1$ under covering relations. The result follows.
  \end{proof}
\end{Corollary}

\medskip

\begin{figure}
  \includegraphics[height=15cm]{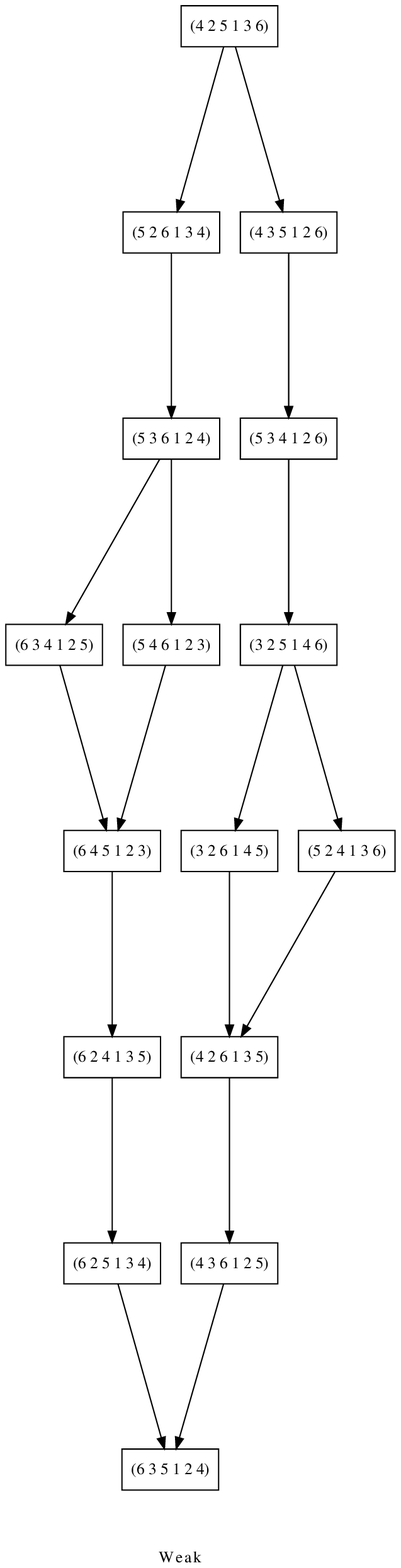}
  \hspace{2cm}
    \includegraphics[height=15cm]{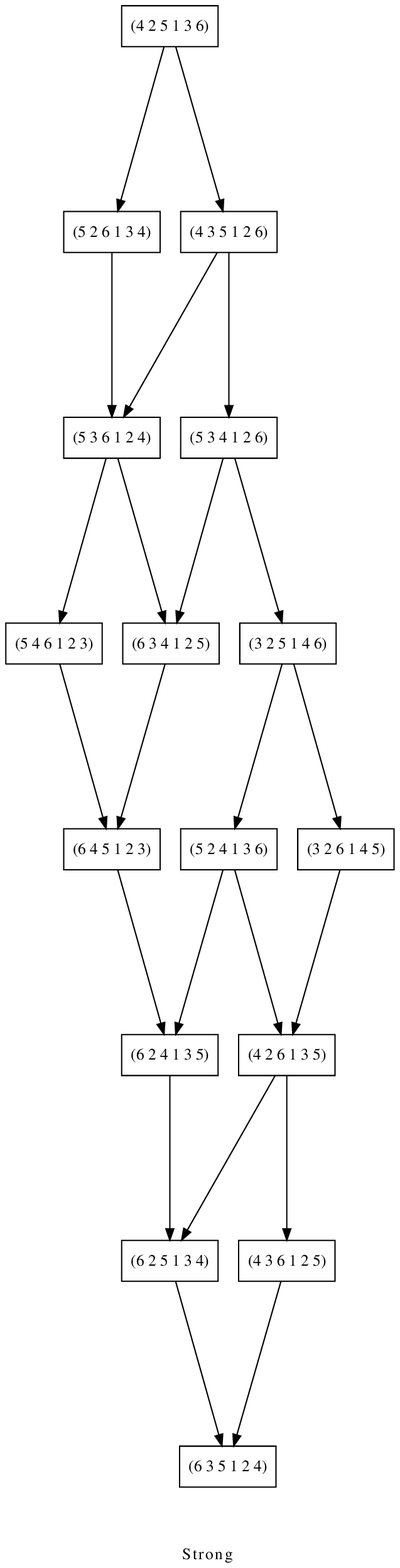}
    \caption{Hasse diagram of the Weak SYT Poset and the Strong SYT
      Poset of $\lambda=(3,2,1)$. Each tableau is represented by its
      row reading word in these pictures. }
  \label{fig:weak.order}
\end{figure}

In \Cref{fig:weak.order}, we show an example of both the Weak SYT
Poset and the Strong SYT poset for $\lambda=(3,2,1)$.  More examples
of these partial orders are given at
\url{https://sites.math.washington.edu/~billey/papers/syt.posets}.

\begin{Remark}
  Observe that both the positive and negative rotation rules apply
  equally well to any skew shape tableaux in $\SYT(\lambda/\nu)$.  The
  block rules apply to skew shape tableaux as well when $T_{z}$ is a
    straight shape tableau.  However, in order to define the analogous
    posets on $\SYT(\lambda/\mu)$, one must include additional block
    moves.  This is part of an ongoing project.
  \end{Remark}

\begin{Remark}
Lascoux--Sch\"utzenberger \cite{LS7} defined an operation called
\textit{cyclage} on semi\-standard tableaux, which decreases
\textit{cocharge} by $1$. The \textit{cyclage poset} on the set of
semi\-standard tableaux arises from applying cyclage in all possible
ways. Cyclage preserves the \textit{content}, i.e.~the number of
$1$'s, $2$'s, etc.  See also \cite[Sect. 4.2]{Shimozono-Weyman.2000}.
Restricting to standard tableaux, cocharge coincides with $\maj$, so
the cyclage poset on $\SYT(n)$ is ranked by $\maj$. However, cyclage
does not necessarily preserve the shape, so it does not suffice to
prove \Cref{thm:zeros}.  For example, restricting the cyclage poset to
$\SYT(32)$ gives a poset which has two connected components and is
not ranked by $\maj$, while both of our poset structures on
$\SYT(32)$ are chains.  A reviewer of \cite{BKS.FPSAC.2020}
posed an interesting question: is
there any relation between the cyclage poset covering relations
restricted to $\SYT(\lambda)$ and the two ranked poset structures used
to prove \Cref{thm:zeros}?  We have not found one, but such a
connection would be interesting if found.
\end{Remark}

\section{Internal zeros for $\des$ on $\SYT(\lambda)$}\label{sec:AER}
The results of \Cref{sec:internal_zeros} show that
$\SYT(\lambda)^{\maj}(q)$ almost never has internal zeros.
Adin--Elizalde--Roichman analogously considered the internal zeros of
the descent number generating functions $\SYT(\lambda/\nu)^{\des}(q)$
where $\des(T )$ is the number of descents in a tableau $T$.

\begin{Question}{{\cite[Problem~7.5]{1801.00044}}}\label{que:des_zeros}
  Is $\{\des(T) : T \in \SYT(\lambda/\nu)\}$ an interval consisting of
  consecutive integers, for any skew shape $\lambda/\nu$? That is,
  does $\SYT(\lambda/\nu)^{\des}(q)$ ever have internal zeros?
\end{Question}

The minimum and maximum descent numbers are easily described
as fol\-lows. The argument involves constructions similar to the sequences of
vertical and horizontal strips used in \Cref{def.min.max.tabs}.

\begin{Lemma}{{\cite[Lemma~3.7(1)]{1801.00044}}}
  Let $\lambda/\nu$ be a skew shape with $n$ cells.  Let $c$ be the maximum length of a
  column and $r$ be the maximum length of a row. Then
  \begin{align*}
    \min \des(\SYT(\lambda/\nu))
      &= c-1 \\
    \max \des(\SYT(\lambda/\nu))
      &= (n-1) - (r-1) = n-r.
  \end{align*}
\end{Lemma}

Indeed, it is easy to see that $\minmaj(\lambda)$ constructed as in
\Cref{def.min.max.tabs} has $\lambda'_1 - 1 = c-1$ descents, and
symmetrically that $\maxmaj(\lambda)$ has $n-r$ descents.
The arguments in \Cref{sec:internal_zeros} consequently resolve
\Cref{que:des_zeros} affirmatively in the straight-shape case.

\begin{Corollary}\label{cor:AER}
  For $\lambda \vdash n$, we have
    \[ \{\des(T) : T \in \SYT(\lambda)\} = \{\lambda_1'-1, \lambda_1', \ldots,
          n - \lambda_1 - 1, n - \lambda_1\}. \]
  In particular, $\SYT(\lambda)^{\des}(q)$ has no internal zeros.

  \begin{proof}
    First suppose $\lambda$ is not a rectangle with at least two rows and columns.
    Iterating the $\varphi$ map
    creates a chain from $\minmaj(\lambda)$ to $\maxmaj(\lambda)$.
    At each step, $\varphi$ either applies a rotation rule or a block rule.
    Rotation rules preserve descent number. Block rules always
    increase the descent number by exactly $1$. Since $\minmaj(\lambda)$
    and $\maxmaj(\lambda)$ have the minimum and maximum number of descents
    possible, the result follows.

    If $\lambda$ is a rectangle with at least two rows and columns, it is easy
    to see that the unique tableau of major index $b(\lambda)+2$ has
    exactly one more descent than $\minmaj(\lambda)$. The result
    follows as before by iterating the $\varphi$ map.
  \end{proof}
\end{Corollary}

\begin{Remark}
  The same argument shows that $\SYT(\lambda)^{\maj - \des}(q)$
  also has no internal zeros. Indeed, applying a rotation rule increases $\maj - \des$
  by $1$ while fixing $\des$, and applying a B-rule fixes $\maj - \des$ and increases
  $\des$ by $1$. In this sense, the strong or weak posets $P(\lambda)$ and $Q(\lambda)$
  have a $\bZ \times \bZ$ ranking given by $(\maj - \des, \des)$.
\end{Remark}

\section{Internal zeros for fake degrees of $C_m \wr  S_n$}\label{ssec:zeros_wreath}
In this section, we classify which irreducible representations appear
in which degrees of the corresponding coinvariant algebras for all
finite groups of the form $C_m \wr S_n$.  The goal is to classify when
the fake degrees $b_{\underline{\lambda}, k} \neq 0$.  We will use the
following helpful lemma which is straightforward to prove.

\begin{Lemma}\label{lem:products} Suppose that $f$ and $g$ are polynomials in
  $\mathbb{Z}[q]$ with non-negative coefficients, that $f$ has no
  internal zeros and has at least two non-zero co\-ef\-fi\-cients, and that
  $g$ has no adjacent internal zeros. Then, $fg$ has no internal
  zeros.
\end{Lemma}

\begin{Lemma}\label{lem:block_internal_zeros}
  Let $\underline{\lambda}=(\lambda^{(1)},\ldots, \lambda^{(m)})$ be a
  sequence of partitions.  The poly\-no\-mi\-al
  $\SYT(\underline{\lambda})^{\maj}(q)$ has no internal zeros except
  when $\underline{\lambda}$ has a single non-empty block $\lambda^{(i)}$
  which is a rectangle with at least two rows and columns. In this
  latter case, the only internal zero up to symmetry occurs at
  $k=b(\lambda^{(i)}) + 1$.

  \begin{proof}
    If $\underline{\lambda}$ has only one nonempty partition, then the
    characterization of internal zeros follows from \Cref{thm:zeros},
    so assume $\underline{\lambda}$ has two or more nonempty
    partitions.
    From \Cref{thm:diag_maj}, we have
    \begin{equation}\label{eq:gmd.case.zeros}
     \SYT(\underline{\lambda})^{\maj}(q) =
      \binom{n}{|\lambda^{(1)}|, \ldots, |\lambda^{(m)}|}_q
      \prod_{i=1}^m \SYT(\lambda^{(i)})^{\maj}(q).
      \end{equation}
      By MacMahon's \Cref{thm:macmahon}, we observe that the
      $q$-multinomial coefficients have no internal zeros.
      Furthermore, the  $q$-multinomial in
      \eqref{eq:gmd.case.zeros} is not constant whenever
      $\underline{\lambda}$ has two or more non-empty partitions.
      From \Cref{thm:zeros}, we know $\SYT(\lambda^{(i)})^{\maj}(q)$
      has no \textit{adjacent} internal zeros for any
      $1\leq i \leq m$.  Consequently, by \Cref{lem:products}, the
      overall product in \eqref{eq:gmd.case.zeros} has no internal
      zeros.
  \end{proof}
\end{Lemma}

\begin{Theorem}\label{thm:block_internal_zeros}
  Let $\underline{\lambda}$ be a sequence of $m$ partitions with
  $|\underline{\lambda}| = n$, and assume
  $g^{\underline{\lambda}}(q) = \sum_k b_{\underline{\lambda}, k}
  q^k$.  Then for $k \in \mathbb{Z}$,
  $b_{\underline{\lambda}, k} \neq 0$ if and only if
    \[ \frac{k - b(\alpha(\underline{\lambda}))}{m} - b(\underline{\lambda})
        \in \left\{0, 1, \ldots, \binom{n+1}{2} - \sum_{c \in \underline{\lambda}} h_c\right\}
             \setminus \cD_{\underline{\lambda}}, \]
  where $\cD_{\underline{\lambda}}$ is empty unless $\underline{\lambda}$ has
  a single non-empty partition $\lambda^{(i)}$ which is a rectangle with at least two rows
  and columns, in which case
    \[ \cD_{\underline{\lambda}}
        = \left\{1, \binom{n+1}{2} - \sum_{c \in \lambda^{(i)}} h_c - 1\right\}. \]

  \begin{proof}
    By \Cref{thm:block_exponents},
    \[ g^{\underline{\lambda}}(q) = q^{b(\alpha(\underline{\lambda}))}
      \SYT(\underline{\lambda})^{\maj}(q^m) \] which implies
    $b_{\underline{\lambda}, k} \neq 0$ only if
    $k-b(\alpha(\underline{\lambda}))$ is a multiple of $m$.  By
    \Cref{lem:block_internal_zeros}, we know
    $\SYT(\underline{\lambda})^{\maj}(q)$ has either no internal zeros
    or internal zeros at degree $1+b(\underline{\lambda})$ and degree
    one less than the maximal major
    index for $\underline{\lambda}$ in the case of a rectangle with at
    least 2 rows and columns.  By \eqref{eq:min.maj} and
    \eqref{eq:max.maj}, the minimal major index for
    $\underline{\lambda}$ is $ b(\underline{\lambda}):= \sum_i b(\lambda^{(i)})$, and the maximal
    major index is
    $b(\underline{\lambda}) + \binom{|\underline{\lambda}|+1}{2} -
    \sum_{c \in \underline{\lambda}} h_c$.  Hence, the result follows.
  \end{proof}
\end{Theorem}

\begin{Corollary}\label{cor:internal_zeros_B}
  In type $B_n$, the irreducible representation indexed by
  $(\lambda, \mu)$ with $|\lambda| + |\mu| = n$ appears in degree $k$
  of the coinvariant algebra of $G(2, 1, n)$ if and only if
    \[ \frac{k - |\mu|}{2} - b(\lambda) - b(\mu) \in \left\{0, 1, \ldots,
        \binom{n+1}{2} - \sum_{c \in \lambda} h_c - \sum_{c' \in \mu} h_{c'}\right\}
        \setminus \cD_{(\lambda, \mu)}. \]
\end{Corollary}

\begin{Example}\label{ex:B6}
  Consider the type $B_6$ case, where $m=2,d=1,n=6$.  For
  $\underline{\lambda}=((2), (3 1))$, we have
  $b(\underline{\lambda})=1$ and
  $$g^{\underline{\lambda}}(q)=
  q^{26} + 2q^{24} + 4q^{22} + 5q^{20} + 7q^{18} + 7q^{16} + 7q^{14} +
  5q^{12} + 4q^{10} + 2q^{8} + q^{6}.$$
For
  $\underline{\mu}=(\varnothing, (3 3))$, we have $b(\underline{\mu})=3$
  and
  $$g^{\underline{\mu}}(q)=
q^{24} + q^{20} + q^{18} + q^{16} + q^{12}.
$$
In both cases, the nonzero coefficients are determined by
\Cref{cor:internal_zeros_B}.
\end{Example}

\section{Deformed Gaussian Multinomial Coefficients}\label{sec:deformed}
We now turn our attention to extending
\Cref{thm:block_internal_zeros} to general Shephard--Todd
groups $G(m, d, n)$.  We begin by introducing a deformation of the
$q$-multinomial coefficients arising from
\Cref{thm:block_exponents.2} in the special case when
$\underline{\lambda}=((\alpha_1), (\alpha_2), \ldots, (\alpha_m))$
is a sequence of one row partitions.  After several lemmas, we give an
alternative formulation for these deformed $q$-multinomials in terms
of inversion generating functions on words with a bounded first
letter.

\begin{Definition}\label{def:q.binomial.deformation}
  Let $\alpha = (\alpha_1, \ldots, \alpha_m) \vDash n$ be a weak
  composition of $n$ with $m$ parts.  Recall the long cycle
  $\sigma_m = (1,2,\ldots,m) \in S_m$, so
    \[ \sigma_m \cdot \alpha = (\alpha_m, \alpha_1, \alpha_2, \ldots, \alpha_{m-1}). \]
  Let $d \mid m$, $\tau = \sigma_m^{m/d}$, and
  $C_d =\langle \tau \rangle = \langle \sigma_m^{m/d} \rangle$ so
  $C_d$ acts on length $m$
  compositions by $(m/d)$-fold cyclic rotations as in
  \Cref{def:cononical.tabs}.  Set
  \begin{align}
    \bbinom{n}{\alpha}_{q;d}
      \coloneqq \frac{\sum_{\sigma \in C_d} q^{b(\sigma \cdot \alpha)}}{[d]_{q^{nm/d}}}
      \binom{n}{\alpha}_{q^m}
  \end{align}
  where
  \[ b(\alpha) \coloneqq \sum_{i=1}^m (i-1) \alpha_i. \]
 Note that when $q=1$, we have  \ $\bbinom{n}{\alpha}_{1;d} = \binom{n}{\alpha}$,
 and when $d=1$, we have
  $\bbinom{n}{\alpha}_{q;1} = q^{b(\alpha)} \binom{n}{\alpha}_{q^m}$, where $m$ is
  the number of parts of $\alpha$. As usual, we also
  write $\bbinom{n}{k}_{q; d} \coloneqq \bbinom{n}{k, n-k}_{q; d}
  = \bbinom{n}{n-k, k}_{q; d}$, where $m=2$ in this case.
  Note that $\bbinom{n}{\alpha}_{q; d}$ is invariant under the $C_d$-action
  on $\alpha$, though this is not typically true of general permutations
  of $\alpha$.
\end{Definition}

\begin{figure}[ht]
    \includegraphics[height=4cm]{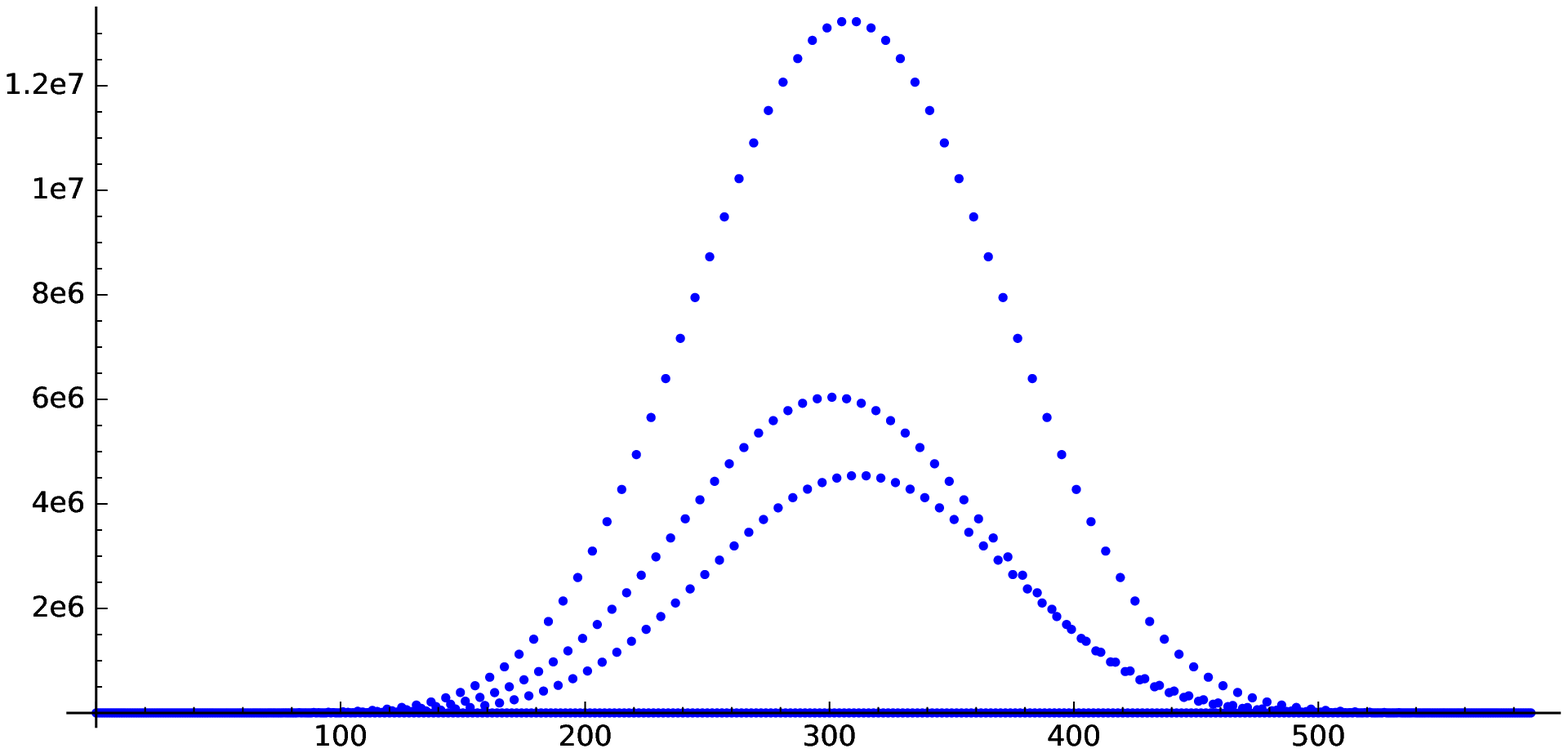}
  \caption{A plot of the coefficients for the deformed $q$-multinomial
    $\bbinom{n}{\alpha}_{q; d}$ with $\alpha= (2, 1, 3, 1, 4, 5)$ and $d=3$.}
  \label{fig:213145.3}
\end{figure}

\begin{Example}
  Observe that $\binom{n}{\alpha}_{q^m}$ alone is generally not
  divisible by $[d]_{q^{nm/d}}$.  For example, if $n = 5$, $\alpha = (2,1,1,1)$,
  and $d=2$, we have
  $$\binom{5}{2,1,1,1}_{q^4}= q^{36} + 3q^{32} + 6q^{28} + 9q^{24} +
  11q^{20} + 11q^{16} + 9q^{12} + 6q^{8} + 3q^{4} + 1$$
  which is not divisible by $[2]_{q^{5\cdot4/2}}=q^{10} + 1$. However,
  $\sum_{\sigma \in C_d} q^{b(\sigma \cdot \alpha)}=q^8+q^{6}$ and
  $(q^8+q^{6})\binom{5}{2,1,1,1}_{q^4}$ is divisible by $q^{10}+1$
  giving
  \begin{align*}\bbinom{5}{2,1,1,1}_{q;2}&= q^{34} + q^{32} + 3q^{30} + 3q^{28} +
  6q^{26} + 5q^{24} + 8q^{22} \\&+ 6q^{20} + 8q^{18} + 5q^{16} + 6q^{14} +
    3q^{12} + 3q^{10} + q^{8} + q^{6} .
  \end{align*}
  See \Cref{fig:213145.3} for a larger example.
\end{Example}

\begin{Lemma}\label{lem:r_diff}
  Given $\alpha = (\alpha_1, \ldots, \alpha_m)\vDash n$, we have

    \[ b(\sigma_m \cdot \alpha) - b(\alpha) = n - m\alpha_m, \]
and
    \[ b(\tau \cdot \alpha) - b(\alpha)
        = nm/d - m( \alpha_m+
        \alpha_{m-1}+ \cdots + \alpha_{m-m/d+1}). \]

  \begin{proof}
    The second claim follows by iterating the first for
    $\tau = \sigma_m^{m/d}$. For the first, we have
    \begin{align*}
      b(\sigma_m \cdot \alpha) - b(\alpha)
        &= (\alpha_1 + 2\alpha_2 + \cdots + (m-1)\alpha_{m-1}) \\
        &\quad
             - (\alpha_2 + 2\alpha_3 + \cdots + (m-1)\alpha_m) \\
        &= \alpha_1 + \alpha_2 + \cdots + \alpha_{m-1} - (m-1)\alpha_m,
    \end{align*}
    which simplifies to $n - m\alpha_m$.
  \end{proof}
\end{Lemma}

If $\alpha \vDash n$, let $\dec{i}\alpha$ be the vector obtained
from $\alpha$ by decreasing $\alpha_i$ by $1$.  Extend the
definition of $\binom{n}{\alpha}_q$ to $m$-tuples of integers by
declaring $\binom{n}{\alpha}_q=0$ if any $\alpha_i$ is negative.  The
following lemma is well known but we include a proof for
completeness.

\begin{Lemma}\label{lem:q_mult_sum}
  We have the following recurrence for $q$-multinomial coefficients,
  \[ \binom{n}{\alpha_1, \ldots, \alpha_m}_q = \sum_{i=1}^m
    q^{\alpha_1 + \cdots + \alpha_{i-1}} \binom{n-1}{\dec{i}\alpha}_q.
    \]

  \begin{proof} By MacMahon's Theorem, the left-hand side is the
    inversion number generating function on length $n$ words with
    $\alpha_i$ copies of the letter $i$ for each $i$. If the first
    letter in such a word is $i$, the number of inversions involving
    the first letter is $\alpha_1 + \alpha_2 + \cdots + \alpha_{i-1}$,
    from which the result quickly follows.
  \end{proof}
\end{Lemma}

The non-trivial deformation of the $q$-binomial coefficients in
\Cref{def:q.binomial.deformation} has the following more explicit form.
In particular, these rational functions are
polynomials with non-negative integer coefficients that satisfy a
Pascal-type formula.

\begin{Lemma}\label{lem:q_binom_def}
  In the case $d=m=2$, we have
\begin{equation}\label{eq:q.pascal.def} \bbinom{n}{k}_{q; 2} = \frac{q^k + q^{n-k}}{1 + q^n} \binom{n}{k}_{q^2}
        = q^{n-k} \binom{n-1}{k-1}_{q^2} + q^k \binom{n-1}{k}_{q^2}
        \in \bZ_{\geq 0}[q].
      \end{equation}
  \begin{proof}
    The first equality is immediate from
    \Cref{def:q.binomial.deformation}. For the second, we use the
    well-known ``$q$-Pascal'' identities
    \[ \binom{n}{k}_q = q^k \binom{n-1}{k}_q + \binom{n-1}{k-1}_q =
      \binom{n-1}{k}_q + q^{n-k} \binom{n-1}{k-1}_q, \]
    which arise from \Cref{lem:q_mult_sum}.  Thus,
    \[q^k \binom{n}{k}_{q^2} = q^k \binom{n-1}{k}_{q^2} + q^{n+n-k}
      \binom{n-1}{k-1}_{q^2}
      \]
      and
      \[
        q^{n-k} \binom{n}{k}_{q^2} = q^{n+k} \binom{n-1}{k}_{q^2} + q^{n-k}\binom{n-1}{k-1}_{q^2}.
      \]
      Hence,
      \[(q^k + q^{n-k}) \binom{n}{k}_{q^2} = (1+q^n)\left(q^k
          \binom{n-1}{k}_{q^2} + q^{n-k}\binom{n-1}{k-1}_{q^2}\right)
      \]
      so the second equality in \eqref{eq:q.pascal.def} holds.
  \end{proof}
\end{Lemma}

We next generalize \Cref{lem:q_binom_def} to all
$\bbinom{n}{\alpha}_{q; d}$ for any $\alpha \vDash n$.  The proof that
follows is independent of \Cref{thm:block_exponents.2}, which
can also be used to prove they are polynomials with non-negative integer coefficients.

\begin{Theorem}\label{thm:q_d_mult_sum}
  Let $\alpha$ be a weak composition of $n$ with $m$ parts, and let $d \mid m$. Then
  \[ \bbinom{n}{\alpha}_{q; d} = \sum_{\sigma \in C_d} q^{b(\sigma
      \cdot \alpha)} \sum_{v=1}^{m/d} q^{m \cdot ((\sigma \cdot
      \alpha)_1 + \cdots + (\sigma \cdot \alpha)_{v-1})}
    \binom{n-1}{\dec{v}(\sigma \cdot \alpha)}_{q^m}. \]
  In particular, $\bbinom{n}{\alpha}_{q; d}$ is a polynomial with
  non-negative coefficients.

  \begin{proof}
    Observe from the definition that
    $\binom{n}{\alpha}_{q}= \binom{n}{\sigma \cdot \alpha}_{q}$ for
    any $\sigma \in C_d$.  Thus, by \Cref{lem:q_mult_sum}, we can
    rewrite the numerator of $\bbinom{n}{\alpha}_{q;d}$ as
    \begin{align*}
      \sum_{\sigma \in C_d} q^{b(\sigma \cdot \alpha)} \binom{n}{\alpha}_{q^m}
        &= \sum_{j=1}^d q^{b(\tau^j \cdot \alpha)} \binom{n}{\tau^j \cdot \alpha}_{q^m} \\
        &= \sum_{j=1}^d \sum_{i=1}^m
             q^{\epsilon(i, j, \alpha)} \binom{n-1}{\dec{i}(\tau^j \cdot \alpha)}_{q^m}
    \end{align*}
    where
    \begin{equation}\label{eq:epsilon}
       \epsilon(i, j, \alpha) \coloneqq
          b(\tau^j \cdot \alpha) + m \cdot ((\tau^j \cdot \alpha)_1 + \cdots +
          (\tau^j \cdot \alpha)_{i-1}).
    \end{equation}
    It is straightforward to check that $\dec{i}(\sigma_m \cdot \alpha)
    = \sigma_m \cdot \dec{i-1}\alpha$, so that $\dec{i}(\tau^j \cdot \alpha)
    = \tau^j \cdot \dec{i-jm/d}\alpha$, where indices are taken modulo $m$.
    Thus,
    \[ \sum_{\sigma \in C_d} q^{b(\sigma \cdot \alpha)}
      \binom{n}{\alpha}_{q^m} = \sum_{i=1}^m \sum_{j=1}^d
      q^{\epsilon(i, j, \alpha)}
      \binom{n-1}{\dec{i-jm/d}\alpha}_{q^m}. \]
    Group the terms on
    the right according to the value $i-jm/d \equiv_m t \in [m]$.
    Note that $j\in [d]$ could be equivalently represented as
    $j \in \bZ/d$, though $i \in [m]$ cannot be treated similarly here.
    One may check that the set of
    $(i, j) \in [m] \times \bZ/d$ such that
    $i-jm/d \equiv_m t$ can be described as
    \[ \{(t+sm/d, s) : s \in [-p_t, d-1-p_t]\} \]
    where $t = p_t m/d + v_t$ for some unique $p_t \in [0,d-1]$ and
    $v_t \in [m/d]$.
    Consequently,
      \[ \sum_{\sigma \in C_d} q^{b(\sigma \cdot \alpha)} \binom{n}{\alpha}_{q^m}
          = \sum_{t=1}^m \left(\sum_{s=-p_t}^{d-1-p_t} q^{\epsilon(t+sm/d, s, \alpha)}\right)
          \binom{n-1}{\dec{t}\alpha}_{q^m}. \]

        Next, we evaluate the incremental change
        \[ \epsilon(t+(s+1)m/d, s+1, \alpha) - \epsilon(t+sm/d, s,
        \alpha) \] for a given $s$.  Let $\beta = \tau^s \cdot
        \alpha$. By \Cref{lem:r_diff},
    \begin{align*}
      b(\tau^{s+1} \cdot \alpha) - b(\tau^s \cdot \alpha)
        &= b(\tau \cdot \beta) - b(\beta) \\
        &= nm/d - m \cdot ( \beta_m + \cdots + \beta_{m-m/d+1}).
    \end{align*}
    We also find
    \begin{align*}
      (\tau \cdot \beta)_1 + & \cdots + (\tau \cdot
                               \beta)_{t+(s+1)m/d-1} = \beta_{m-m/d+1} +\cdots + \beta_m+\beta_{1}+ \cdots + \beta_{t+sm/d-1}
    \end{align*}
so
        \begin{align*}
      (\tau \cdot \beta)_1 + & \cdots + (\tau \cdot \beta)_{t+(s+1)m/d-1}
      - \beta_1 - \cdots - \beta_{t+sm/d-1} \\
           &= \beta_{m-m/d+1} + \cdots + \beta_m.
    \end{align*}
    Combining these observations,
    \begin{align*}
      \epsilon(t+(s+1)m/d, s+1, \alpha) - \epsilon(t+sm/d, s, \alpha)
        = nm/d.
    \end{align*}
    It follows that
    \begin{align*}
      \sum_{s=-p_t}^{d-1-p_t} q^{\epsilon(t+sm/d, s, \alpha)}
        &= q^{\epsilon(t-p_tm/d, -p_t, \alpha)} [d]_{q^{nm/d}} \\
        &= q^{\epsilon(v_t, -p_t, \alpha)} [d]_{q^{nm/d}}.
    \end{align*}
    Since we have a bijection $[0, d-1] \times [m/d] \to [m]$ given by
    $(p, v) \mapsto pm/d + v$, we have
    \begin{equation}\label{eq:pv}
      \sum_{\sigma \in C_d} q^{b(\sigma \cdot \alpha)} \binom{n}{\alpha}_{q^m}
        = [d]_{q^{nm/d}} \sum_{p=0}^{d-1} \sum_{v=1}^{m/d} q^{\epsilon(v, -p, \alpha)}
              \binom{n-1}{\dec{v+pm/d}\alpha}_{q^m},
    \end{equation}
    proving the polynomiality of $\bbinom{n}{\alpha}_{q; d}$.

    We can further refine \eqref{eq:pv}.   From \eqref{eq:epsilon}, we
    observe that
      \[ \epsilon(v, -p, \alpha) = \epsilon(v, 0, \tau^{-p} \cdot
    \alpha), \]
    and since $\tau=\sigma_m^{m/d}$,  we have
     \[\dec{v+pm/d}\alpha
         = \tau^{p} \cdot \dec{v}(\tau^{-p} \cdot \alpha).\]
So,
      \[ q^{\epsilon(v, -p, \alpha)}
              \binom{n-1}{\dec{v+pm/d}\alpha}_{q^m}  =
              q^{\epsilon(v, 0, \tau^{-p} \cdot \alpha)} \binom{n-1}{\dec{v}(\tau^{-p} \cdot \alpha)}_{q^m}, \]
which implies
    \begin{align*}
      \sum_{\sigma \in C_d} q^{b(\sigma \cdot \alpha)} \binom{n}{\alpha}_{q^m}
      &= [d]_{q^{nm/d}} \sum_{\sigma \in C_d} q^{b(\sigma \cdot \alpha)}
        \sum_{v=1}^{m/d} q^{m \cdot ((\sigma \cdot \alpha)_1
        + \cdots + (\sigma \cdot \alpha)_{v-1})}
        \binom{n-1}{\dec{v}(\sigma \cdot \alpha)}_{q^m}.
    \end{align*}
    The result follows by dividing by $[d]_{q^{nm/d}}$.
  \end{proof}
\end{Theorem}

In light of \Cref{thm:q_d_mult_sum}, we define the following
polynomials.

\begin{Definition}\label{def:pam}
  Let $\alpha=(\alpha_1,\ldots, \alpha_m) \vDash n$, and say
  $1 \leq k \leq m$.  Define the $\alpha,k$-\emph{partial sum multinomial}
  by
\[ p_{\alpha}^{(k)}(q) = \sum_{i=1}^k q^{\alpha_1 + \cdots +
    \alpha_{i-1}} \binom{n-1}{\dec{i}\alpha}_q.
\]
\end{Definition}

\begin{Remark}
  By \Cref{lem:q_mult_sum}, $p_{\alpha}^{(m)}=\binom{n}{\alpha}_q$, and more generally
  the same argument shows that
  \begin{equation}\label{eq:pak.words}
    p_{\alpha}^{(k)}(q) = \{w \in \W_\alpha : w_1 \leq k\}^{\inv}(q)
  \end{equation}
  is an inversion number generating function.
\end{Remark}

It is very well-known that the multinomial coefficients can be written as a product of
binomial coefficients. More generally, $q$-multinomial coefficients can be written
as a product of $q$-binomial coefficients. This holds true even for the
$\alpha, k$-partial sum multinomials as follows.

\begin{Lemma}\label{lem:pam_prod}
  Let $\alpha=(\alpha_1, \ldots, \alpha_m) \vDash n$ and $1 \leq k \leq m$.
  We have
  \begin{equation}\label{eq:pam_prod}
    p_{\alpha}^{(k)}(q) = \prod_{i=1}^k
      \binom{\alpha_1 + \cdots + \alpha_i}{\alpha_i}_q \cdot \prod_{i=k+1}^m
      \binom{\alpha_1 + \cdots + \alpha_i - 1}{\alpha_i}_q.
  \end{equation}
\end{Lemma}

\begin{proof}
  Recall that
  $p_{\alpha}^{(k)}(q) = \{w \in \W_\alpha : w_1 \leq k\}^{\inv}(q)$.
  Partition the set $\{w \in \W_\alpha : w_1 \leq k\}$ into
  $\binom{n-1}{\alpha_1+\cdots+\alpha_k-1, \alpha_{k+1}, \ldots, \alpha_m}$
  subsets according to the placement of all $k+1,k+2,\ldots,m$'s in positions
  $2,3,\ldots,n$.  For each such placement, there are
  $\binom{\alpha_1+\cdots+\alpha_k}{\alpha_1,\ldots,\alpha_k}$ ways to
  place numbers $1,2,\ldots, k$ in the remaining
  positions. Since each inversion in a word $w \in \W_\alpha$ is
  between two letters $\leq k$, between two letters $\geq k+1$, or between a
  letter $\leq k$ and a letter $\geq k+1$, it follows that
     \begin{equation}\label{eq:pak.prod}
       \{w \in \W_\alpha : w_1 \leq
       k\}^{\inv}(q) =       \binom{\alpha_1 + \cdots +
         \alpha_k}{\alpha_1,\ldots, \alpha_k}_q
       \binom{n-1}{\alpha_1+\cdots+ \alpha_k-1,\alpha_{k+1},\ldots, \alpha_m}_q
     \end{equation}
     by MacMahon's \Cref{thm:macmahon}.  Factoring each
     $q$-multinomial in \eqref{eq:pak.prod} into $q$-binomials gives
     \eqref{eq:pam_prod}.
  \end{proof}

\begin{Corollary}\label{cor:pam_unimodal}
  Let $\alpha=(\alpha_1, \ldots, \alpha_m) \vDash n$ and
  $1 \leq k \leq m$.  The $\alpha, k$-partial sum multinomial
  $p_\alpha^{(k)}(q)$ is symmetric and unimodal.

  \begin{proof}
    A result of Andrews \cite[Thm.~3.9]{MR0557013} states that the
    product of symmetric, unimodal polynomials with non-negative
    coefficients is symmetric and unimodal with non-negative
    coefficients. The $q$-binomials are symmetric with non-negative
    coefficients, and it is a well-known, non-trivial fact that they
    are also unimodal. See \cite{Zeilberger.1986} for a combinatorial
    proof of this fact and further historical references. The result now follows from
    \Cref{lem:pam_prod}.
  \end{proof}
\end{Corollary}

\begin{Lemma}\label{lem:pam_constant}
  Let $\alpha=(\alpha_1, \ldots, \alpha_m) \vDash n$ and
  $1 \leq k \leq m$.  Then $p_\alpha^{(k)}(q) \neq 0$ if and only if
  $\alpha_1 + \cdots + \alpha_k > 0$.  In this case,
  $p_\alpha^{(k)}(q)$ has constant coefficient $1$, degree
  $D_\alpha - \alpha_{k+1} - \cdots - \alpha_m$ where
  $D_\alpha = \binom{n}{2} - \sum \binom{\alpha_i}{2}$ is the degree
  of $\binom{n}{\alpha}_q$, and has no internal zeros. Furthermore,
  $p_{\alpha}^{(k)}(q)$ is non-constant except when
  \begin{itemize}
  \item $\alpha_1 + \cdots + \alpha_k = 0$, in which case $p_{\alpha}^{(k)}(q)=0$;
  \item $\alpha_1 + \cdots + \alpha_k = 1$ and
    $\alpha_i=n-1$ for some $i > k$, in which case $p_{\alpha}^{(k)}(q)=1$; or
  \item $\alpha_i=n$ for some $i \leq k$, in which case $p_{\alpha}^{(k)}(q)=1$.
  \end{itemize}

  \begin{proof}
    Each claim follows easily from the fact that
    $p_{\alpha}^{(k)}(q)$ is the inversion generating function for
    $\{w \in \W_\alpha: w_1\leq k\}$. Alternatively, one may use
    \Cref{lem:pam_prod}.
  \end{proof}
\end{Lemma}

We have the following summary statement.

\begin{Corollary}\label{cor:deformed.multinomial.rec}
  Let $\alpha$ be weak composition of $n$ with $m$ parts, and let
  $d \mid m$. Let
  $\underline{\alpha}=((\alpha_1),(\alpha_2), \ldots,(\alpha_m))
  $ be the corresponding sequence of one row partitions.  Then,
  \[
    \bbinom{n}{\alpha}_{q;d} =
      \frac{\sum_{\sigma \in C_d} q^{b(\sigma \cdot \alpha)}}{[d]_{q^{nm/d}}}
      \binom{n}{\alpha}_{q^m} \ = \
    \sum_{\sigma \in C_d} q^{b(\sigma \cdot \alpha)} p_{\sigma \cdot
      \alpha}^{(m/d)}(q^m) \ = \
    \frac{d}{\#\{\underline{\alpha}\}^d}\
    \gmdn{\alpha}(q).
  \]
\end{Corollary}

\begin{proof}
  The first equality is just the definition.  The second equality
  follows from \Cref{thm:q_d_mult_sum}.  The third equality follows
  from \Cref{thm:block_exponents.2} and the fact that
  $$\sum_{\sigma \in C_d} q^{b(\sigma \cdot \alpha)} =
  \frac{d}{\#\{\underline{\alpha}\}^d}\ \left(\{\underline{\alpha}\}^d\right)^{b\circ \alpha}(q).$$
  \end{proof}

  \begin{Remark}
    We note that since
    $ \bbinom{n}{\alpha}_{q;d} = \frac{d}{\#\{\underline{\alpha}\}^d}\
    \gmdn{\alpha}(q)$, we knew from Stembridge's work that the
    deformed multinomial coefficients are polynomials in $q$ even
    though they are defined as rational functions.  Our proof in
    The\-o\-rem \ref{thm:q_d_mult_sum} gives an alternate, direct proof of this fact without
    going through representation theory.  Furthermore, we use the
    summation formula in Cor\-ol\-lary \ref{cor:deformed.multinomial.rec} to
    characterize the internal zeros of $\gmdn{\lambda}(q)$ in the next section.
\end{Remark}

\section{Internal zeros for $G(m, d, n)$}\label{sec:other_types}

We can now extend the results of \Cref{ssec:zeros_wreath} to all
Shephard--Todd groups $G(m, d, n)$.  We thus give a remarkably simple
and completely general de\-scrip\-tion for which irreducible
representations appear in which degrees of the coinvariant algebras of
essentially arbitrary complex reflection groups.  Recall the notation
established in \Cref{sub:complex.reflection.groups}.  Let
$\{\underline{\lambda}\}^d$ be the orbit of $\underline{\lambda}$ under
$(m/d)$-fold rotations in $C_d =\langle \sigma_m^{m/d} \rangle$.

\begin{Definition}
  Given a sequence $\underline{\lambda} = (\lambda^{(1)}, \ldots, \lambda^{(m)})$
  with $|\underline{\lambda}| = n$, let
    \[ \alpha(\underline{\lambda}) \coloneqq
        (|\lambda^{(1)}|, \ldots, |\lambda^{(m)}|) \vDash n. \]
  Similarly, let $\underline{\alpha}(\underline{\lambda})$ be the length $m$ sequence
  of partitions whose $i$th partition
  is the single row partition of size $|\lambda^{(i)}|$.
\end{Definition}

The map $\underline{\alpha}$ may not be injective on
$\{\underline{\lambda}\}^d$, though it has constant fiber sizes since
$\underline{\alpha}$ is $C_d$-equivariant.  For example, when
$m=4, d=4$, we have
\begin{align*}
  \underline{\alpha} \colon
       ((2), \varnothing, (1^2), \varnothing) &\mapsto ((2), \varnothing, (2), \varnothing) \\
       (\varnothing, (2), \varnothing, (1^2)) &\mapsto (\varnothing, (2), \varnothing, (2)) \\
       ((1^2), \varnothing, (2), \varnothing) &\mapsto ((2), \varnothing, (2), \varnothing) \\
       (\varnothing, (1^2), \varnothing, (2)) &\mapsto (\varnothing, (2), \varnothing, (2)).
\end{align*}

Generalizing \Cref{thm:diag_maj}, we have the following corollary of
Stembridge's \Cref{thm:block_exponents.2} and
\Cref{def:q.binomial.deformation}.

\begin{Corollary}\label{lem:g_la_d}
  Let $\underline{\lambda}$ be a sequence of $m$ partitions with
  $|\underline{\lambda}| = n$. Let $\{\underline{\lambda}\}^d$ be the orbit of
  $\underline{\lambda}$ under $(m/d)$-fold cyclic rotations. Then
  \begin{align*}
    \gmdn{\lambda}(q)
    &=
      \frac{\#\{\underline{\lambda}\}^d}{d} \cdot
      \bbinom{n}{\alpha(\underline{\lambda})}_{q; d} \cdot
      \prod_{i=1}^m \SYT(\lambda^{(i)})^{\maj}(q^m).
  \end{align*}

  \begin{proof}
    By \Cref{thm:block_exponents.2},
    \begin{align*}
      \gmdn{\lambda}(q)
        &= \frac{\left(\{\underline{\lambda}\}^d\right)^{b\circ \alpha}(q)}{[d]_{q^{nm/d}}}
          \SYT(\underline{\lambda})^{\maj}(q^m).
    \end{align*}
    We have
    $$\left(\{\underline{\lambda}\}^d\right)^{b\circ \alpha}(q)=
    \frac{\#\{\underline{\lambda}\}^d}{\#\{\underline{\alpha}(\underline{\lambda})\}^d}
    \cdot \left(\{\underline{\alpha}(\underline{\lambda})\}^d\right)^{b\circ \alpha}(q)$$
    and
    $$
    \left(\{\underline{\alpha}(\underline{\lambda})\}^d\right)^{b\circ \alpha}(q)=
    \frac{\#\{\underline{\alpha}(\underline{\lambda})\}^d}{d}
    \sum_{\sigma \in C_d} q^{b(\sigma \cdot
      \alpha(\underline{\lambda}))} .
$$
Consequently, using  \Cref{thm:diag_maj} and \Cref{def:q.binomial.deformation}, we have
        \begin{align*}
      \gmdn{\lambda}(q)
        &= \frac{\#\{\underline{\lambda}\}^d}{d} \cdot
             \frac{\sum_{\sigma \in C_d} q^{b(\sigma \cdot
      \alpha(\underline{\lambda}))}}{[d]_{q^{nm/d}}} \cdot
             \binom{n}{\alpha(\underline{\lambda})}_{q^m} \cdot
             \prod_{i=1}^m \SYT(\lambda^{(i)})^{\maj}(q^m) \\
        &= \frac{\#\{\underline{\lambda}\}^d}{d} \cdot
                \bbinom{n}{\alpha(\underline{\lambda})}_{q; d} \cdot
          \prod_{i=1}^m
             \SYT(\lambda^{(i)})^{\maj}(q^m).
    \end{align*}
  \end{proof}
\end{Corollary}

We will now prove the general classification theorem for nonzero fake
degrees as mentioned in the introduction.  The reader may find it
useful to compare the statement to the type $A$ case in
\Cref{thm:zeros} and the $C_m \wr S_n$ case in
\Cref{thm:block_internal_zeros}.

\begin{Theorem}\label{thm:block_internal_zeroes.2}
  Let $\underline{\lambda}$ be a sequence of $m$ partitions with
  $|\underline{\lambda}| = n \geq 1$, let $d \mid m$, and let
  $\{\underline{\lambda}\}^d$ be the orbit of
  $\underline{\lambda}$ under the group $C_d$ of $(m/d)$-fold cyclic rotations. Then
  $b_{\{\underline{\lambda}\}^d, k} \neq 0$ if and only if for some $\underline{\mu}
  \in \{\underline{\lambda}\}$ we have $|\mu^{(1)}| + \cdots + |\mu^{(m/d)}| > 0$
  and
    \[ \frac{k - b(\alpha(\underline{\mu}))}{m} - b(\underline{\mu})
        \in \left\{0, 1, \ldots, |\mu^{(1)}| + \cdots + |\mu^{(m/d)}|
        + \binom{n}{2} - \sum_{c \in \underline{\mu}} h_c\right\}
        \setminus \cD_{\underline{\mu}; d}. \]
  Here   $\cD_{\underline{\mu}; d}$ is empty unless either
  \begin{enumerate}[(1)]
    \item $\underline{\mu}$ has a partition $\mu$ of size $n$; or
    \item $\underline{\mu}$ has a partition $\mu$ of size $n-1$ and
      $|\mu^{(1)}|  + \cdots + |\mu^{(m/d)}| = 1$,
  \end{enumerate}
  where in both cases $\mu$ must be a rectangle with at least two rows and
  columns. In case (1), we have
    \[ \cD_{\underline{\mu}; d}
        \coloneqq \left\{1, \binom{n+1}{2} - \sum_{c \in \mu} h_c -1 \right\}, \]
  and in case (2) we have
    \[ \cD_{\underline{\mu}; d}
        \coloneqq \left\{1, \binom{n}{2} - \sum_{c \in \mu} h_c \right\}. \]

\end{Theorem}

  \begin{proof} Let
    $\alpha=\alpha(\underline{\lambda})$.
    Using \Cref{lem:g_la_d} and \Cref{cor:deformed.multinomial.rec}, we have
    \begin{equation}\label{eq:gmdn.master}
      \gmdn{\lambda}(q)
        = \frac{\#\{\underline{\lambda}\}^d}{d} \cdot \sum_{\sigma \in C_d} q^{b(\sigma
          \cdot \alpha)}
          p_{\sigma \cdot \alpha}^{(m/d)}(q^m) \cdot
          \prod_{i=1}^m \SYT(\lambda^{(i)})^{\maj}(q^m).
    \end{equation}
    Thus, we consider the locations of the nonzero terms in
    \begin{equation}\label{eq:block_internal_zeros.1}
      p_{\sigma \cdot \alpha}^{(m/d)}(q)
      \prod_{i=1}^m \SYT(\lambda^{(i)})^{\maj}(q).
    \end{equation}
    Recall that $p_{\sigma \cdot \alpha}^{(m/d)}(q)=0$ whenever
    $(\sigma \cdot \alpha)_1 + \cdots + (\sigma \cdot \alpha)_{m/d}=0$, so assume
    $(\sigma \cdot \alpha)_1 + \cdots + (\sigma \cdot \alpha)_{m/d}>0$.  Since
    $\SYT(\lambda^{(i)})^{\maj}(q) \neq 0$ for all partitions
    $\lambda^{(i)}$, we can also assume
    \eqref{eq:block_internal_zeros.1} is not zero.

    By \Cref{lem:pam_constant}, $p_{\sigma \cdot \alpha}^{(m/d)}(q)$
    has no internal zeros, degree
    $ \binom{n}{2} -\sum \binom{\alpha_i}{2} -(\sigma \cdot \alpha)_{m/d+1}
    - \cdots - (\sigma \cdot \alpha)_m$, and constant term
    $p_{\sigma \cdot \alpha}^{(m/d)}(0)=1$.  Thus the minimal degree
    term of \eqref{eq:block_internal_zeros.1} with nonzero coefficient
    is $q^{b(\underline{\lambda})}$ by \eqref{eq:min.maj}, and the
    maximal degree term is $q$ to the power
    \begin{equation}\label{eq:max-deg.g}
      \binom{n}{2} -\sum \binom{\alpha_i}{2} -(\sigma \cdot \alpha)_{m/d+1} -
    \cdots - (\sigma \cdot \alpha)_m
      + \sum \mathrm{deg}(\SYT(\lambda^{(i)})^{\maj}(q)).
    \end{equation}
    Since $\alpha_i=|\lambda^{(i)}|$, we know by \eqref{eq:max.maj}
    that
    $\mathrm{deg}(\SYT(\lambda^{(i)})^{\maj}(q)) =
    \binom{\alpha_i}{2}
    - b(\lambda^{(i)\, \prime})$, so \eqref{eq:max-deg.g} simplifies to
    \begin{equation}\label{eq:max-deg.2}
      \binom{n}{2}  -(\sigma \cdot \alpha)_{m/d+1} -
    \cdots - (\sigma \cdot \alpha)_m
      - b(\underline{\lambda}')
    \end{equation}
    where
    $b(\underline{\lambda}') :=\sum_i b(\lambda^{(i)\, \prime})$.
    From \eqref{eq:max.maj}, we also know
    $\binom{n}{2}  - b(\underline{\lambda}') = b(\underline{\lambda})
    +  \binom{n+1}{2} -  \sum_{c \in \underline{\lambda}} h_c$
and $\sigma \cdot \alpha \vDash n$,  so we
conclude that the maximal degree of \eqref{eq:block_internal_zeros.1} is
    \begin{equation}\label{eq:max-deg.3}
b(\underline{\lambda})+ (\sigma \cdot \alpha)_{1} +\cdots + (\sigma \cdot \alpha)_{m/d} +
      \binom{n}{2} -  \sum_{c \in \underline{\lambda}} h_c.
    \end{equation}

    If $p_{\sigma \cdot \alpha}^{(m/d)}(q)\neq 1$, then the product in
    \eqref{eq:block_internal_zeros.1} also has no internal zeros by
    \Cref{thm:zeros} and \Cref{lem:products}.
    The cases where $p_{\sigma \cdot \alpha}^{(m/d)}(q)=1$ are listed
    in \Cref{lem:pam_constant}.  By assumption,
    $(\sigma \cdot \alpha)_1 + \cdots + (\sigma \cdot \alpha)_{m/d}>0$ so the
    remaining cases are when
    $(\sigma \cdot \alpha)_1 + \cdots + (\sigma \cdot \alpha)_{m/d}=1$ and some
    $\alpha_i=n-1$ for $i > m/d$, or $\alpha_i=n$ for some $i \leq m/d$. In either of the
    remaining cases,  determining the nonzero coefficients of
    \eqref{eq:block_internal_zeros.1} reduces to the case of a single
    partition described in \Cref{thm:zeros}.

    To finish the proof, we observe from
    \eqref{eq:gmdn.master} that
    $b_{\{\underline{\lambda}\}^d, k} \neq 0$ if and only if there
    exists some
    $\sigma \in C_d$ such that
    $(\sigma \cdot \alpha)_1 + \cdots + (\sigma \cdot \alpha)_{m/d} > 0$ and the corresponding
    product
$$q^{b(\sigma \cdot
  \alpha)} p_{\sigma \cdot \alpha}^{(m/d)}(q^m) \cdot \prod_{i=1}^m
\SYT(\lambda^{(i)})^{\maj}(q^m)$$ has nonzero coefficient of $q^k$.
Thus, by our analysis of the location of nonzero coefficients in
\eqref{eq:block_internal_zeros.1}, we observe that
$b_{\{\underline{\lambda}\}^d, k} \neq 0$ if and only if
    \[ \frac{k - b(\sigma \cdot \alpha)}{m} -b(\underline{\lambda})
      \in \left\{0,1,
        \ldots,
 (\sigma \cdot \alpha)_{1} +\cdots + (\sigma \cdot \alpha)_{m/d} +
      \binom{n}{2} -  \sum_{c \in \underline{\lambda}} h_c
      \right\}
      \setminus \cD_{\underline{\mu}; d}
       \]
       where $\underline{\mu}=\sigma \cdot \underline{\lambda}$.  Observing that
       $b(\underline{\lambda}) =b(\underline{\mu})$,\
       $\sum_{c \in \underline{\lambda}} h_c = \sum_{c \in
         \underline{\mu}} h_c$, and $b(\sigma \cdot \alpha)
         = b(\alpha(\underline{\mu}))$ completes the proof of the theorem.
  \end{proof}

\begin{Corollary}\label{cor:internal_zeros_D}
  In type $D_n$, an irreducible with orbit $\{\lambda, \mu\}$ where $|\lambda| + |\mu| = n$
  appears in degree $k$ of the coinvariant algebra of $G(2, 2, n)$ if and only if
  either $\lambda \neq \varnothing$ and
    \[ \frac{k-|\mu|}{2} - b(\lambda) - b(\mu)
        \in \left\{0, 1, \ldots, |\lambda| + \binom{n}{2}
        - \sum_{c \in \lambda} h_c - \sum_{c' \in \mu} h_{c'}\right\}
        \setminus \cD_{(\lambda, \mu); 2} \]
  or $\mu \neq \varnothing$ and
    \[ \frac{k-|\lambda|}{2} - b(\lambda) - b(\mu)
        \in \left\{0, 1, \ldots, |\mu| + \binom{n}{2}
        - \sum_{c \in \lambda} h_c - \sum_{c' \in \mu} h_{c'}\right\}
        \setminus \cD_{(\mu, \lambda); 2}. \]
\end{Corollary}

\begin{Example}
  Consider the type $D_6$ case, where $m=2,d=2,n=6$.  For
  $\underline{\lambda}=((2), (3 1))$, we have
  $b(\underline{\lambda})=1$ and
  $$g^{\{\underline{\lambda}\}^2}(q)=
q^{20} + 3q^{18} + 6q^{16} + 8q^{14} + 9q^{12} + 8q^{10} + 6q^{8} + 3q^{6} + q^{4} .$$ On the other hand, if
  $\underline{\mu}=(\varnothing, (3 3))$ then $b(\underline{\mu})=3$
  and
  $$g^{\{\underline{\mu}\}^2}(q)=q^{18} + q^{14} + q^{12} + q^{10} + q^{6} . $$
  In both cases, the nonzero coefficients are determined by
  \Cref{cor:internal_zeros_D}.  One may notice that
  $g^{\{\underline{\mu}\}^2}(q) \neq g^{\{\underline{\mu}\}^1}(q)$,
  which appeared in \Cref{ex:B6}.  However, for
  $\underline{\nu}=( (3 3),\varnothing)$, one can check
  $\{\underline{\mu}\}^2 = \{\underline{\nu}\}^2$ and
$g^{\{\underline{\mu}\}^2}(q) = g^{\{\underline{\nu}\}^2}(q)=
g^{\{\underline{\nu}\}^1}(q)$.
\end{Example}


\section{Future Work}\label{sec:future}

A sequence $a_0, a_1, a_2, \ldots$ is \textit{parity-unimodal} if
$a_0, a_2, a_4, \ldots$ and $a_1, a_3, a_5, \ldots$ are each
unimodal. Stucky \cite[Thm.~1.3]{stucky2018cyclic} recently showed that the
$q$-Catalan polynomials, namely $\SYT((n, n))^{maj}(q)$ up to a
$q$-shift, are parity-unimodal. The argument involves constructing an
$\fsl_2$-action on rational Che\-rednik algebras. See [\S3.1, Haiman94]
for a prototype of the argument in a highly related context.
Recent work of Gaetz--Gao \cite{MR4042823} constructed an $\fsl_2$-module on
$\bC S_n$ and strongly related work of Hamaker--Pechenik--Speyer--Weigandt
\cite{ALCO_2020__3_2_301_0}
constructed an $\fsl_2$-module on $R_n$, though neither of these structures
are capable of producing internal zeros and they do not
respect the isotypic decomposition. Nonetheless, based on
Stucky's result, our internal zeros classification, and a brute-force
check for $n \leq 50$, we conjecture the following.

\begin{Conjecture} The fake-degree polynomials $f^\lambda(q)$ are parity-unimodal for all $\lambda$.
\end{Conjecture}

When $W$ is a Weyl group, the Hilbert series of the coinvariant
algebra $R_W$ is symmetric and unimodal by the Hard Lefschetz Theorem
since $R_W$ presents the cohomology of the associated flag variety
$G/B$. One referee asked the following interesting question.

\begin{quote}
Is there an algebraic or geometric witness to the fact that the
$f^\lambda(q)$ so rarely have internal zeros in their coefficient
sequences?  More precisely, let $\ell = c_1 x_1 + \cdots + c_n x_n$ be
a linear form in $\bC[x_1, \ldots, x_n]$ with $c_i \neq c_j$ whenever
$i \neq j$. These are precisely the Lefschetz elements when $W =
S_n$. If $V_\lambda = \oplus_{d=0}^{n(n-1)/2} (V_\lambda)_d$ is the
$\lambda$-isotypical component of $R_n$, decomposed by polynomial
degree, for each $d$ we have a linear map \[ (V_\lambda)_d \too{\times
\ell} (R_n)_{d+1} \too{\epsilon_\lambda} (V_\lambda)_{d+1}, \] where
we first multiply by $\ell$ and then we act by the Young symmetrizer
$\epsilon_\lambda \in \bC[S_n]$. Initial computer verifications
suggest this composite linear map is nonzero whenever $f^\lambda(q)$
does not have an internal zero, at least for some special choices of
the coefficients $c_i$. Could this correspond to some property in
ge\-om\-e\-try?
\end{quote}

\begin{Remark}
  We note that Stanley used the Hard Lefschetz Theorem to prove that
  Bruhat orders are rank symmetric, rank unimodal, and have a symmetric
  chain decomposition, so they are Sperner \cite{MR804700}.  This
  theorem was part of our motivation for defining the Weak SYT Poset.
  We were looking for a subposet of the Strong SYT Poset which is a
  symmetric chain decomposition, though we did not find
  one. The Weak SYT Poset is not a disjoint union of chains in general, though
  it is the most natural subposet we could find. The
  names ``Weak'' SYT Poset and ``Strong'' SYT Poset simply imply that one is a
  subposet of the other in the same way that the weak order is a
  subposet of ``strong'' Bruhat order.
\end{Remark}

\section*{Acknowledgments}\label{sec:ack}

We would like to thank Christian Stump, William McGovern,
Andrew Oha\-na, Alejandro Morales, Greta Panova, Mihael Perman,
Martin Rai\v{c}, Victor Reiner,
Sheila Sundaram, Vasu Tewari, Lauren Williams and Alex Woo for
helpful discussions related to this work. We also thank both anonymous
referees for their very careful and helpful comments and one of the
referees for suggesting the Hard Lefschetz Theorem approach
in \Cref{sec:future}.

\bibliography{refs}{}

\newcommand{\etalchar}[1]{$^{#1}$}
\begin{thebibliography}{HPSW20}

\bibitem[AER18]{1801.00044}
Ron Adin, Serge Elizalde, and Yuval Roichman.
\newblock Cyclic descents for near-hook and two-row shapes, 2018.
\newblock arXiv:1801.00044.

\bibitem[And76]{MR0557013}
George~E. Andrews.
\newblock {\em The {T}heory of {P}artitions}.
\newblock Addison-Wesley Publishing Co., Reading, Mass.-London-Amsterdam, 1976.
\newblock Encyclopedia of Mathematics and its Applications, Vol. 2.

\bibitem[AS17]{Ahlbach-Swanson.2017}
Connor Ahlbach and Joshua~P. Swanson.
\newblock Refined cyclic sieving.
\newblock {\em S\'em. Lothar. Combin.}, 78B:Art. 48, 12, 2017.

\bibitem[BKS20a]{bks2}
Sara~C. Billey, Matja\v{z} Konvalinka, and Joshua~P. Swanson.
\newblock Asymptotic normality of the major index on standard tableaux.
\newblock {\em Adv. in Appl. Math.}, 113:101972, 36, 2020.
\newblock Online: arXiv:1905.00975.

\bibitem[BKS20b]{BKS.FPSAC.2020}
Sara~C. Billey, Matja\v{z} Konvalinka, and Joshua~P. Swanson.
\newblock On the distribution of the major index on standard {Y}oung tableaux.
\newblock {\em S\'em. Lothar. Combin.}, 2020.
\newblock Extended abstract accepted as a talk for FPSAC 2020. To appear.

\bibitem[Car75]{carlitz.1975}
L.~Carlitz.
\newblock A combinatorial property of {$q$}-{E}ulerian numbers.
\newblock {\em Amer. Math. Monthly}, 82:51--54, 1975.

\bibitem[Car89]{Carter.89}
Roger~W. Carter.
\newblock {\em Simple {G}roups of {L}ie {T}ype}.
\newblock Wiley Classics Library. John Wiley \& Sons, Inc., New York, 1989.
\newblock Reprint of the 1972 original, A Wiley-Interscience Publication.

\bibitem[Cli37]{Clifford.1937}
Alfred~H. Clifford.
\newblock Representations induced in an invariant subgroup.
\newblock {\em Ann. of Math. (2)}, 38(3):533--550, 1937.

\bibitem[ER15]{Ehrenborg.Readdy.2015}
Richard Ehrenborg and Margaret Readdy.
\newblock A poset view of the major index.
\newblock {\em Adv. in Appl. Math.}, 62:1--14, 2015.

\bibitem[Foa68]{Foata}
D.~Foata.
\newblock {On the Netto inversion number of a sequence}.
\newblock {\em Proc.\ Amer.\ Math.\ Soc.}, 19:236--2479--1130, 1968.

\bibitem[FRT54]{Frame-Robinson-Thrall.1954}
J.~S. Frame, G.~de~B. Robinson, and R.~M. Thrall.
\newblock The hook graphs of the symmetric groups.
\newblock {\em Canadian J. Math.}, 6:316--324, 1954.

\bibitem[GG20]{MR4042823}
Christian Gaetz and Yibo Gao.
\newblock A combinatorial {$\mathfrak{sl}_2$}-action and the {S}perner property
  for the weak order.
\newblock {\em Proc. Amer. Math. Soc.}, 148(1):1--7, 2020.

\bibitem[Gre55]{MR0072878}
J.~A. Green.
\newblock The characters of the finite general linear groups.
\newblock {\em Trans. Amer. Math. Soc.}, 80:402--447, 1955.

\bibitem[HHL05]{HHL2005}
J.~Haglund, M.~Haiman, and N.~Loehr.
\newblock A combinatorial formula for {M}acdonald polynomials.
\newblock {\em J. Amer. Math. Soc.}, 18(3):735--761, 2005.

\bibitem[HPSW20]{ALCO_2020__3_2_301_0}
Zachary Hamaker, Oliver Pechenik, David~E Speyer, and Anna Weigandt.
\newblock Derivatives of {S}chubert polynomials and proof of a determinant
  conjecture of {S}tanley.
\newblock {\em Algebraic Combinatorics}, 3(2):301--307, 2020.

\bibitem[Hum90]{Hum}
James~E. Humphreys.
\newblock {\em Reflection {G}roups and {C}oxeter {G}roups}, volume~29 of {\em
  Cambridge Studies in Advanced Mathematics}.
\newblock Cambridge University Press, Cambridge, 1990.

\bibitem[Kly74]{klyachko74}
A.~A. Klyachko.
\newblock Lie elements in the tensor algebra.
\newblock {\em Siberian Mathematical Journal}, 15(6):914--920, 1974.

\bibitem[Knu73]{Knuth}
D.~E. Knuth.
\newblock {\em {The Art of Computer Programming}}, volume~3.
\newblock Addison--Wesley, Reading, MA, 1973.

\bibitem[LS81]{LS7}
Alain Lascoux and Marcel-P. Sch{\"u}tzenberger.
\newblock Le mono\"\i de plaxique.
\newblock In {\em Noncommutative structures in algebra and geometric
  combinatorics ({N}aples, 1978)}, volume 109 of {\em Quad. ``Ricerca Sci.''},
  pages 129--156. CNR, Rome, 1981.

\bibitem[Lus77]{Lusztig.77}
George Lusztig.
\newblock Irreducible representations of finite classical groups.
\newblock {\em Invent. Math.}, 43(2):125--175, 1977.

\bibitem[Mac13]{MacMahon.1913}
P.~A. MacMahon.
\newblock The indices of permutations and the derivation therefrom of functions
  of a single variable associated with the permutations of any assemblage of
  objects.
\newblock {\em Amer. J. Math.}, 35(3):281--322, 1913.

\bibitem[Mac17]{MR1576566}
P.~A. MacMahon.
\newblock Two applications of general theorems in combinatory analysis: (1) to
  the theory of inversions of permutations; (2) to the ascertainment of the
  numbers of terms in the development of a determinant which has amongst its
  elements an arbitrary number of zeros.
\newblock {\em Proc. London Math. Soc.}, S2-15(1):314, 1917.

\bibitem[MPP15]{1512.08348}
Alejandro Morales, Igor Pak, and Greta Panova.
\newblock Hook formulas for skew shapes {I}. $q$-analogues and bijections,
  2015.
\newblock arXiv:1512.08348.

\bibitem[{OEI}18]{oeis}
{OEIS Foundation Inc.}
\newblock The {O}n-{L}ine {E}ncyclopedia of {I}nteger {S}equences, 2018.
\newblock Online. \url{http://oeis.org}.

\bibitem[RS{\etalchar{+}}18]{findstat}
Martin Rubey, Christian Stump, et~al.
\newblock {FindStat} - {T}he combinatorial statistics database.
\newblock \url{http://www.FindStat.org}, 2018.

\bibitem[RSW04]{Reiner-Stanton-White.CSP}
V.~Reiner, D.~Stanton, and D.~White.
\newblock The cyclic sieving phenomenon.
\newblock {\em J. Combin. Theory Ser. A}, 108(1):17--50, 2004.

\bibitem[Sag91]{Sagan.1991}
B.~Sagan.
\newblock {\em {The Symmetric Group}}.
\newblock Wadsworth, Inc., Belmont, CA, 1991.

\bibitem[Spe35]{Specht.1935}
W.~Specht.
\newblock Die irreduziblen darstellungen der symmetrischen gruppe.
\newblock {\em Mathematische Zeitschrift}, 39:696--711, 1935.

\bibitem[ST54]{Shephard-Todd.54}
G.~C. Shephard and J.~A. Todd.
\newblock Finite unitary reflection groups.
\newblock {\em Canadian J. Math.}, 6:274--304, 1954.

\bibitem[Sta79]{stanley.1979}
Richard~P. Stanley.
\newblock Invariants of finite groups and their applications to combinatorics.
\newblock {\em Bull. Amer. Math. Soc. (N.S.)}, 1(3):475--511, 1979.

\bibitem[Sta80]{Stanley.1980}
Richard~P. Stanley.
\newblock Weyl groups, the hard {Lefschetz} theorem, and the {Sperner}
  property.
\newblock {\em SIAM J. Alg. Disc. Meth.}, 1(2):168--184, 1980.

\bibitem[Sta84]{MR804700}
Richard~P. Stanley.
\newblock Combinatorial applications of the hard {L}efschetz theorem.
\newblock In {\em Proceedings of the {I}nternational {C}ongress of
  {M}athematicians, volume\ 1, 2 ({W}arsaw, 1983)}, pages 447--453. PWN,
  Warsaw, 1984.

\bibitem[Sta99]{ec2}
R.~P. Stanley.
\newblock {\em Enumerative {C}ombinatorics. volume 2}, volume~62 of {\em
  Cambridge Studies in Advanced Mathematics}.
\newblock Cambridge University Press, Cambridge, 1999.

\bibitem[Sta12]{ec1}
Richard~P. Stanley.
\newblock {\em Enumerative {C}ombinatorics. volume 1}, volume~49 of {\em
  Cambridge Studies in Advanced Mathematics}.
\newblock Cambridge University Press, Cambridge, second edition, 2012.

\bibitem[Ste51]{steinberg.1951}
R.~Steinberg.
\newblock A geometric approach to the representations of the full linear group
  over a {G}alois field.
\newblock {\em Trans. Amer. Math. Soc.}, 71:274--282, 1951.

\bibitem[Ste89]{stembridge89}
John~R. Stembridge.
\newblock On the eigenvalues of representations of reflection groups and wreath
  products.
\newblock {\em Pacific J. Math.}, 140(2):353--396, 1989.

\bibitem[Stu18]{stucky2018cyclic}
Eric Stucky.
\newblock Cyclic sieving, necklaces, and bracelets.
\newblock arXiv:1812.04578, 2018.

\bibitem[SW00]{Shimozono-Weyman.2000}
Mark Shimozono and Jerzy Weyman.
\newblock Graded characters of modules supported in the closure of a nilpotent
  conjugacy class.
\newblock {\em European J. Combin.}, 21(2):257--288, 2000.

\bibitem[SW10]{Shareshian-Wachs.2010}
John Shareshian and Michelle~L. Wachs.
\newblock Eulerian quasisymmetric functions.
\newblock {\em Adv. Math.}, 225(6):2921--2966, 2010.

\bibitem[Swa18]{s17}
Joshua~P. Swanson.
\newblock On the existence of tableaux with given modular major index.
\newblock {\em Algebraic Combinatorics}, 1(1):3--21, 2018.

\bibitem[You77]{Young.collected.work}
Alfred Young.
\newblock {\em The {C}ollected {P}apers of {A}lfred {Y}oung (1873--1940)}.
\newblock University of Toronto Press, Toronto, Ont., Buffalo, N. Y., 1977.

\bibitem[Zei89]{Zeilberger.1986}
Doron Zeilberger.
\newblock Kathy {O}'{H}ara's constructive proof of the unimodality of the
  {G}aussian polynomials.
\newblock {\em Amer. Math. Monthly}, 96(7):590--602, 1989.

\end{thebibliography}
\bibliographystyle{alpha}

\end{document}